\documentclass[12pt]{article}
\usepackage[dvipsnames]{xcolor}
\usepackage{amsfonts}
\usepackage{setspace}
\usepackage{bm}
\usepackage{bbm}
\usepackage{amssymb}
\usepackage{indentfirst}
\usepackage{geometry}
\usepackage{mathtools}
\usepackage{tocloft}

\DeclarePairedDelimiter\abs{\lvert}{\rvert}

\makeatletter
\let\oldabs\abs
\def\abs{\@ifstar{\oldabs}{\oldabs*}}

\newcommand*{\citena}[1]{%
\begingroup
\color{Green}
\romannumeral-`\x 
\setcitestyle{numbers}%
\citep{#1}
\endgroup
\ignorespacesafterend
}

\newcommand*{\citesup}[1]{%
\begingroup
\color{Green}
\citep{#1}
\endgroup
\ignorespacesafterend
}

\newcommand*{\eqrefe}[1]{%
\begingroup
(\color{BrickRed}
\romannumeral-`\x 
\setcitestyle{numbers}%
\ref{eq:#1}%
\endgroup
)\ignorespacesafterend
}

\newcommand*{\secrefe}[1]{%
\begingroup
(\color{Aquamarine}
\romannumeral-`\x 
\setcitestyle{numbers}%
\ref{#1}%
\endgroup
)\ignorespacesafterend
}

\newtagform{Tags}[\textcolor{BrickRed}]{\color{Black}(}{)}

\usepackage[numbers, sort]{natbib}
\setcitestyle{super}

\newcommand{\ii}{\bm{i}}

\begin{document}
\title{Generalized Harmonic Numbers}
\author{Jose Risomar Sousa}
\date{October 18, 2018}
\maketitle
\usetagform{Tags}

\begin{abstract}
This paper presents new formulae for the harmonic numbers of order $k$, $H_{k}(n)$, and for the partial sums of two Fourier series associated with them, denoted here by $C^z_{k}(n)$ and $S^z_{k}(n)$. I believe this new formula for $H_{k}(n)$ is an improvement over the digamma function, $\psi$, because it is simpler and it stems from Faulhaber's formula, which provides a closed-form for the sum of powers of the first $n$ positive integers. We demonstrate how to create an exact power series for the harmonic numbers, a new integral representation for $\zeta(2k+1)$ and a new generating function for $\zeta(2k+1)$, among many other original results. The approaches and formulae discussed here are entirely different from solutions available in the literature.
\end{abstract}

\tableofcontents

\section{Introduction}
Although formulae for the harmonic numbers have been known for some time, they are not very simple or useful. For example, a formula due to Euler expresses $H(n)$ as $\int_{0}^{1}(1+x+\dotsb+x^{n-1})\,dx$ for integer $n$, but it is frequently dismissed by scholars, who prefer the approximation $H(n)\sim\log(n)+\gamma$.\\

In this paper we find out how to obtain a more natural and elegant formula for the generalized harmonic numbers: 
\begin{equation}
H_{k}(n)=\sum_{j=1}^{n} \frac{1}{j^k} \nonumber
\end{equation}
\indent This new formula has the advantage of being easier to work with. For example, it can be used to obtain the sum of $H(n)/n^2$ over the positive integers relatively easy.\\

\indent We also show how to obtain the partial sums of two Fourier series associated with $H_{k}(n)$, denoted here by $S^z_{k}(n)$ and $C^z_{k}(n)$, which cover some notable particular cases, such as the alternating harmonic numbers, $C^2_{k}(n)$, and the odd alternating harmonic numbers, $S^4_{k}(n)$ ($S^4_{2}(n)$ converges to Catalan's constant). These two functions are given below (for all integer $k\geq 1$ and complex $z$):
\begin{equation} \nonumber
C^z_{k}(n)=\sum_{j=1}^{n} \frac{1}{j^k}\cos{\frac{2\pi j}{z}} \text{, and } S^z_{k}(n)=\sum_{j=1}^{n} \frac{1}{j^k}\sin{\frac{2\pi j}{z}}
\end{equation}
\indent We create general formulae for $C^z_{k}(n)$ and $S^z_{k}(n)$ and find out their limits as $n$ approaches infinity as a function of Riemann's zeta function. (After looking up previous results in the literature, I found that the limits of $C^z_{2k}(n)$ and $S^z_{2k+1}(n)$ are not new, they are a function of the so-called Bernoulli polynomials\citesup{Abra}, though the limits of $C^z_{2k+1}(n)$ and $S^z_{2k}(n)$ are probably new.)\\

So, to begin, let us recall Faulhaber's formula for the sum of the $i$-th powers of the first $n$ positive integers:
\begin{equation} \nonumber
\sum_{k=1}^{n} {k^i}=\sum_{j=0}^{i} \frac{(-1)^j i! B_j n^{i+1-j}}{(i+1-j)!j!}
\end{equation}
where $B_j$ are the Bernoulli numbers\citesup{Conway}.\\

Since odd Bernoulli numbers are always 0, except for $B_1$, we can
simplify the above formula for even and odd powers as follows:
\begin{equation} \label{eq:soma_pot_par}
\sum_{k=1}^{n} {k^{2i}}=\frac{n^{2i}}{2}+\sum_{j=0}^{i} \frac{(2i)! B_{2j} n^{2i+1-2j}}{(2j)!(2i+1-2j)!} 
\end{equation} 

\begin{equation} \label{eq:soma_pot_impar}
\sum_{k=1}^{n} {k^{2i+1}}=\frac{n^{2i+1}}{2}+\sum_{j=0}^{i} \frac{(2i+1)! B_{2j} n^{2i+2-2j}}{(2j)!(2i+2-2j)!}
\end{equation}

\section{Indicator Function $\mathbbm{1}_{k|n}$}
One key component of the method used to solve the generalized harmonic numbers is the indicator function $\mathbbm{1}_{k|n}$, defined as 1 if $k$ divides $n$, and 0 otherwise. This function and its analog (that will appear in the next section) play a key role in the solution that is presented here:
\begin{equation} \nonumber
\mathbbm{1}_{k|n}=\frac{1}{k}\sum_{j=1}^{k}\cos{\frac{2\pi nj}{k}} 
\end{equation}
\indent A closed-form for $\mathbbm{1}_{k|n}$ can be obtained by means of the so-called Lagrange's trigonometric identities:\\
\begin{equation} \label{eq:k_div_n_trans}
\mathbbm{1}_{k|n}=\frac{1}{2k}\frac{\sin{(2\pi n+\frac{\pi n}{k}})}{\sin{\frac{\pi n}{k}}}-\frac{1}{2k} =\frac{1}{2k}\sin{2\pi n}\cot{\frac{\pi n}{k}}+\frac{\cos{2\pi n}-1}{2k}
\end{equation}
\indent We can also create a power series for $\mathbbm{1}_{k|n}$ by expanding the cosine with Taylor series:
\begin{equation}
\mathbbm{1}_{k|n}=\frac{1}{k}\sum_{j=1}^{k}\cos{\frac{2\pi nj}{k}}=1+\frac{1}{k}\sum_{i=1}^{\infty}\frac{(-1)^i}{(2i)!}\left(\frac{2\pi n}{k}\right)^{2i}\sum_{j=1}^{k}j^{2i} \nonumber
\end{equation}
\indent Now, by replacing the summation of $j^{2i}$ over $j$ with Faulhaber's formula, \eqrefe{soma_pot_par}, we obtain:
\begin{equation} \nonumber
\mathbbm{1}_{k|n}=1+\frac{1}{k}\sum_{i=1}^{\infty}\frac{(-1)^i}{(2i)!}\left(\frac{2\pi n}{k}\right)^{2i}\left(\frac{k^{2i}}{2}+\sum_{j=0}^{i}\frac{(2i)! B_{2j} k^{2i+1-2j}}{(2j)!(2i+1-2j)!}\right) \Rightarrow
\end{equation}

\begin{equation} \label{eq:k_div_n_series1}
\mathbbm{1}_{k|n}=\frac{\cos{2\pi n}-1}{2k}+
\sum_{i=0}^{\infty}(-1)^i(2\pi n)^{2i}\sum_{j=0}^{i}\frac{B_{2j} k^{-2j}}{(2j)!(2i+1-2j)!}
\end{equation} 
\indent From \eqrefe{k_div_n_trans} and \eqrefe{k_div_n_series1}, after re-scaling $n$ to $n/2$, we conclude that:
\begin{equation} \label{eq:even_power_series}
\sum_{i=0}^{\infty} (-1)^i (\pi n)^{2i}\sum_{j=0}^{i} \frac{B_{2j}k^{-2j}}{(2i+1-2j)!(2j)!}=
\frac{1}{2k}\cot{\frac{\pi n}{2k}}\sin{\pi n}
\end{equation}

\subsection{The Analog of $\mathbbm{1}_{k|n}$}
Now, just as we created a power series for $\mathbbm{1}_{k|n}$, we have to create one for its analog, which is the summation:
\begin{equation} \nonumber
\frac{1}{k}\sum_{j=1}^{k}\sin{\frac{2\pi nj}{k}} 
\end{equation}
\indent Again, we can find a closed-form for the above summation using Lagrange's trigonometric identities:
\begin{equation} \label{eq:sum_sine3}
\frac{1}{k}\sum_{j=1}^{k}\sin{\frac{2\pi nj}{k}}=-\frac{1}{2k}\frac{\cos{(2\pi n+\frac{\pi n}{k}})}{\sin{\frac{\pi n}{k}}}+\frac{1}{2k}\cot{\frac{\pi n}{k}}=\frac{\sin{2\pi n}}{2k}+\frac{1}{k}\cot{\frac{\pi n}{k}}\sin^2{\pi n}
\end{equation}
\indent As previously, we can obtain a power series for the above by expanding the sine with Taylor series and making use of \eqrefe{soma_pot_impar}:
\begin{equation} \label{eq:sum_sine4}
\frac{1}{k}\sum_{j=1}^{k}\sin{\frac{2\pi nj}{k}}=\frac{\sin{2\pi n}}{2k}+
\sum_{i=0}^{\infty}(-1)^i(2\pi n)^{2i+1}\sum_{j=0}^{i}\frac{B_{2j} k^{-2j}}{(2i+2-2j)!(2j)!}
\end{equation}
\indent From \eqrefe{sum_sine3} and \eqrefe{sum_sine4}, after re-scaling $n$ to $n/2$, it follows that:
\begin{equation} \label{eq:odd_power_series}
\sum_{i=0}^{\infty}(-1)^i(\pi n)^{2i+1}\sum_{j=0}^{i}\frac{B_{2j} k^{-2j}}{(2i+2-2j)!(2j)!}=
\frac{1}{k}\cot{\frac{\pi n}{2k}}\left(\sin{\frac{\pi n}{2}}\right)^2
\end{equation}

\section{Generalized Harmonic Numbers}

\subsection{Formula Rationale} \label{HN_sin_k_pi}
The rationale to build a formula for $H_{k}(n)$ is to use the Taylor series expansion of $\sin{\pi k}$, and exploit the fact that it is 0 for all integer $k$. We refer to the below as initial equation (note the $k$ in the summation is not the same $k$ used as subscript on $H_{k}(n)$):
\begin{equation} \label{eq:start_iden_sin_k_pi}
\sin{\pi k}=0 \Rightarrow  \pi k=\sum_{i=1}^{\infty}\frac{-(-1)^i (\pi k)^{2i+1}}{(2i+1)!} 
\end{equation}
\indent If we divide both sides of \eqrefe{start_iden_sin_k_pi} by $\pi k^2$ we end up with a power series for $1/k$ about $0$ that only holds for integer $k$ (after all $1/k$ is not analytic at $0$).\\

Besides, on the right-hand side of the resulting equation, the exponents of $k$ are positive integers, allowing to apply Faulhaber's formula mentioned in the introduction. By doing so we end up with a convoluted power series that fortunately can be transformed into an integral by means of the closed-form we derived for $\mathbbm{1}_{k|n}$ (or its analog) using Lagrange's identities. That is a high level summary of the reasoning.\\

To not make this paper long, we only give two fully detailed demonstrations based on the initial equation $\sin{\pi k}=0$, and jump straight to the final formulae in a few other cases, before we state a general formula. We also briefly show how the outcomes change with the choice of different initial equations. 

\subsection{Harmonic Number} \label{HN_1_sin_k_pi}
We start by dividing both sides of \eqrefe{start_iden_sin_k_pi} by $\pi k^2$:
\begin{equation} \label{eq:start_iden_sin_k_pi_1}
\frac{1}{k}=\sum_{i=1}^{\infty}\frac{-(-1)^i \pi^{2i}k^{2i-1}}{(2i+1)!}= 
\sum_{i=0}^{\infty}\frac{(-1)^i \pi^{2i+2}k^{2i+1}}{(2i+3)!} 
\end{equation}
\indent Below we take the sum of \eqrefe{start_iden_sin_k_pi_1} over $k$ and use equation \eqrefe{soma_pot_impar}, thus extending the domain of $H_{1}(n)$ ($H(n)$ for short) to the real numbers, in an analytic continuation:
\begin{equation} \nonumber
H(n)=\sum_{k=1}^{n}\frac{1}{k}=\sum_{i=0}^{\infty}\frac{(-1)^i \pi^{2i+2}}{(2i+3)!}\sum_{k=1}^{n}k^{2i+1}=\sum_{i=0}^{\infty}\frac{(-1)^i \pi^{2i+2}}{(2i+3)!}\left(\frac{n^{2i+1}}{2}+\sum_{j=0}^{i} \frac{(2i+1)! B_{2j} n^{2i+2-2j}}{(2i+2-2j)!(2j)!}\right) 
\end{equation}
\begin{equation} \nonumber
H(n)=\sum_{i=0}^{\infty}\frac{(-1)^i\pi^{2i+2}n^{2i+1}}{2(2i+3)!}
+\sum_{i=0}^{\infty}\frac{(-1)^i\pi^{2i+2}n^{2i+2}}{(2i+3)!}\sum_{j=0}^{i}\frac{(2i+1)! B_{2j} n^{-2j}}{(2i+2-2j)!(2j)!} 
\end{equation}
\indent The first sum is straightforward:
\begin{equation} \nonumber
\sum_{i=0}^{\infty}\frac{(-1)^i\pi^{2i+2}n^{2i+1}}{2(2i+3)!}=-\frac{1}{2\pi n^{2}}\sum_{i=1}^{\infty}\frac{(-1)^i(\pi n)^{2i+1}}{(2i+1)!}=\frac{1}{2\pi n^{2}}\left(\pi n-\sin{\pi n}\right) 
\end{equation}
\indent The second sum is an exact power series for $H(n)-1/(2n)$ and can be rewritten as:
\begin{equation} \label{eq:1_over_k}
\sum_{i=0}^{\infty}\left(-\frac{1}{2i+3}+\frac{1}{2i+2}\right)(-1)^i\pi^{2i+2}n^{2i+2}\sum_{j=0}^{i}\frac{B_{2j} n^{-2j}}{(2i+2-2j)!(2j)!}
\end{equation}
\indent The above sums are tricky, but they can be obtained from \eqrefe{odd_power_series}, one of the formulae derived previously. In order to do that, let us replace $(n,k)$ by $(x,n)$ and define a function $f(x,n)$ such that:
\begin{equation} \label{eq:HN_f_odd_sin_k_pi} 
f(x,n)=\sum_{i=0}^{\infty}(-1)^i(\pi x)^{2i+1}\sum_{j=0}^{i}\frac{B_{2j} n^{-2j}}{(2i+2-2j)!(2j)!}=
\frac{1}{n}\cot{\frac{\pi x}{2n}}\left(\sin{\frac{\pi x}{2}}\right)^2
\end{equation}
\indent To build each piece of \eqrefe{1_over_k}), we start from the above $f(x,n)$.\\

\indent For the first sum, we multiply $f(x,n)$ by $-\pi \cdot x/n$ and integrate with respect to $x$ as below:
\begin{equation} \nonumber
-\frac{1}{n}\sum_{i=0}^{\infty}(-1)^i\pi^{2i+2}\left(\int_{0}^{n} x^{2i+2}\,dx\right)\sum_{j=0}^{i}\frac{B_{2j} n^{-2j}}{(2i+2-2j)!(2j)!} =
-\frac{\pi}{n}\int_{0}^{n} x f(x,n)\,dx  
\end{equation}
\begin{equation} \nonumber
-\sum_{i=0}^{\infty}\frac{(-1)^i\pi^{2i+2}n^{2i+2}}{2i+3}\sum_{j=0}^{i}\frac{B_{2j} n^{-2j}}{(2i+2-2j)!(2j)!}=
-\frac{\pi}{n}\int_{0}^{n} x\frac{1}{n}\cot{\frac{\pi x}{2n}}\left(\sin{\frac{\pi x}{2}}\right)^2dx 
\end{equation}
\indent For the second sum, we multiply $f(x,n)$ by $\pi$ and integrate with respect to $x$ as below:
\begin{equation} \nonumber
\sum_{i=0}^{\infty}(-1)^i\pi^{2i+2}\left(\int_{0}^{n} x^{2i+1}\right)\sum_{j=0}^{i}\frac{B_{2j} n^{-2j}}{(2i+2-2j)!(2j)!}\,dx =\pi\int_{0}^{n} f(x,n)\,dx 
\end{equation}
\begin{equation} \nonumber
\sum_{i=0}^{\infty}\frac{(-1)^i\pi^{2i+2}n^{2i+2}}{2i+2}\sum_{j=0}^{i}\frac{B_{2j} n^{-2j}}{(2i+2-2j)!(2j)!}=\pi\int_{0}^{n} \frac{1}{n}\cot{\frac{\pi x}{2n}}\left(\sin{\frac{\pi x}{2}}\right)^2 dx  
\end{equation}\\
\indent Now, by summing up the two resulting integrals, we obtain the below equivalence (which holds for all real $n$, not just integers):
\begin{multline} \nonumber
\sum_{i=0}^{\infty}\frac{(-1)^i\pi^{2i+2}n^{2i+2}}{(2i+3)!}\sum_{j=0}^{i}\frac{(2i+1)! B_{2j} n^{-2j}}{(2i+2-2j)!(2j)!}=\int_{0}^{n} \frac{\pi(n-x)}{n^2}\cot{\frac{\pi x}{2n}}\left(\sin{\frac{\pi x}{2}}\right)^2 dx 
\\=\pi\int_{0}^{1} u\cot{\frac{\pi (1-u)}{2}}\left(\sin{\frac{\pi n(1-u)}{2}}\right)^2\,du
\end{multline}
\noindent where we have used the transformation $u=1-x/n$.\\

\indent Now, by adding up the simple part (disregarding $\sin{\pi n}$ and changing $u$ for $1-u$), we finally arrive at a formula for $H(n)$:
\begin{equation}
\sum_{k=1}^{n}\frac{1}{k}=\frac{1}{2n}+\frac{\pi}{2}\int_{0}^{1}(1-u)\left(1-\cos{\pi n u}\right)\cot{\frac{\pi u}{2}}\,du 
\end{equation}
\indent This formula has a certain resemblance to Faulhaber's formula, especially the term $1/(2n)$ outside of the integral. If we compare this formula with the one below due to Euler\citesup{Knuth}, based on the digamma function, it seems that the former is more natural and tractable than the latter:
\begin{equation} \nonumber
\indent \sum_{k=1}^{n}\frac{1}{k}=\int_{0}^{1}\frac{1-x^n}{1-x}\,dx=\gamma+\psi(n+1) \text{, where }\gamma \text{ is the Euler-Mascheroni constant.} 
\end{equation}
\indent The two functions approach one another very quickly as $n$ grows large.

\subsection{Harmonic Number of Order 2} \label{HN_2_sin_k_pi} 
We divide both sides of \eqrefe{start_iden_sin_k_pi_1} by $k$:
\begin{equation} \nonumber 
\frac{1}{k^2}=\sum_{i=0}^{\infty}\frac{(-1)^i \pi^{2i+2}k^{2i}}{(2i+3)!} 
\end{equation}
\indent We sum the above over $k$, using equation \eqrefe{soma_pot_par} this time, noting that because equation \eqrefe{soma_pot_par} does not work exactly for $i=0$, we need to make a little correction by adding up $-1/2$:
\begin{equation} \nonumber
H_{2}(n)=\sum_{k=1}^{n}\frac{1}{k^2}=\sum_{i=0}^{\infty}\frac{(-1)^i \pi^{2i+2}}{(2i+3)!}\sum_{k=1}^{n}k^{2i}=-\frac{1}{2}\frac{\pi^2}{3!}+\sum_{i=0}^{\infty}\frac{(-1)^i \pi^{2i+2}}{(2i+3)!}\left(\frac{n^{2i}}{2}+\sum_{j=0}^{i} \frac{(2i)! B_{2j} n^{2i+1-2j}}{(2i+1-2j)!(2j)!}\right) 
\end{equation}
\begin{equation} \nonumber
H_{2}(n)=-\frac{1}{2}\frac{\pi^2}{3!}+\frac{1}{2}\sum_{i=0}^{\infty}\frac{(-1)^i\pi^{2i+2}n^{2i}}{(2i+3)!}
+\sum_{i=0}^{\infty}\frac{(-1)^i\pi^{2i+2}n^{2i+1}}{(2i+3)!}\sum_{j=0}^{i}\frac{(2i)! B_{2j} n^{-2j}}{(2i+1-2j)!(2j)!} 
\end{equation}
\indent The 1st sum, again, is straightforward:
\begin{equation} \nonumber
-\frac{1}{2}\frac{\pi^2}{3!}+\frac{1}{2}\sum_{i=0}^{\infty}\frac{(-1)^i\pi^{2i+2}n^{2i}}{(2i+3)!}=-\frac{1}{2\pi n^3}\sum_{i=2}^{\infty}\frac{(-1)^i(\pi n)^{2i+1}}{(2i+1)!}=\frac{1}{2\pi n^3}\left(\pi n-\frac{(\pi n)^3}{3!}-\sin{\pi n}\right) 
\end{equation}
\indent The 2nd sum can be rewritten as:
\begin{equation} \label{eq:1_over_k_sq} \nonumber
\sum_{i=0}^{\infty}\left(\frac{1}{2(2i+3)}-\frac{1}{2i+2}+\frac{1}{2(2i+1)}\right)(-1)^i\pi^{2i+2}n^{2i+1}\sum_{j=0}^{i}\frac{B_{2j} n^{-2j}}{(2i+1-2j)!(2j)!} 
\end{equation}
\indent The above sums are tricky, but they can be derived from \eqrefe{even_power_series}, another one of the formulae derived previously. In order to do that, let us replace $(n,k)$ by $(x,n)$ and define a function $g(x,n)$ such that:
\begin{equation} \label{eq:HN_g_even_sin_k_pi}
g(x,n)=\sum_{i=0}^{\infty}(-1)^i(\pi x)^{2i}\sum_{j=0}^{i}\frac{B_{2j} n^{-2j}}{(2i+1-2j)!(2j)!}=
\frac{1}{2n}\cot{\frac{\pi x}{2n}}\sin{\pi x} 
\end{equation}
\indent To build each piece of \eqrefe{1_over_k_sq} we start from $g(x,n)$.\\

For the 1st sum, we multiply both sides of $g(x,n)$ by $\pi^2/2\cdot x^2/n^2$ and integrate with respect to $x$ as below:
\begin{equation} \nonumber
\frac{1}{2n^2}\sum_{i=0}^{\infty}(-1)^i\pi^{2i+2}\left(\int_{0}^{n} x^{2i+2}\,dx\right)\sum_{j=0}^{i}\frac{B_{2j} n^{-2j}}{(2i+1-2j)!(2j)!}=\frac{\pi^2}{2n^2}\int_{0}^{n} x^2 g(x,n)\,dx 
\end{equation}
\begin{equation} \nonumber
\frac{1}{2}\sum_{i=0}^{\infty}\frac{(-1)^i\pi^{2i+2}n^{2i+1}}{2i+3}\sum_{j=0}^{i}\frac{B_{2j} n^{-2j}}{(2i+1-2j)!(2j)!}=
\frac{\pi^2}{2n^2}\int_{0}^{n} x^2\frac{1}{2n}\cot{\frac{\pi x}{2n}}\sin{\pi x}\,dx 
\end{equation}
\indent For the 2nd sum, we multiply both sides of $g(x,n)$ by $-\pi^2\cdot x/n$ and integrate with respect to $x$ as below:
\begin{equation} \nonumber
-\frac{1}{n}\sum_{i=0}^{\infty}(-1)^i\pi^{2i+2}\left(\int_{0}^{n} x^{2i+1} \,dx\right)\sum_{j=0}^{i}\frac{B_{2j} n^{-2j}}{(2i+1-2j)!(2j)!}=\frac{-\pi^2}{n}\int_{0}^{n} x g(x,n)\,dx 
\end{equation}
\begin{equation} \nonumber
-\sum_{i=0}^{\infty}\frac{(-1)^i\pi^{2i+2}n^{2i+1}}{2i+2}\sum_{j=0}^{i}\frac{B_{2j} n^{-2j}}{(2i+1-2j)!(2j)!}=
\frac{-\pi^2}{n}\int_{0}^{n} x\frac{1}{2n}\cot{\frac{\pi x}{2n}}\sin{\pi x}\,dx 
\end{equation}
\indent For the 3rd sum, we multiply both sides of $g(x,n)$ by $\pi^2/2$ and integrate with respect to $x$ as below:
\begin{equation} \nonumber
\frac{1}{2}\sum_{i=0}^{\infty}(-1)^i\pi^{2i+2}\left(\int_{0}^{n} x^{2i}\, dx\right)\sum_{j=0}^{i}\frac{B_{2j} n^{-2j}}{(2i+1-2j)!(2j)!}=\frac{\pi^2}{2}\int_{0}^{n} g(x,n)\,dx 
\end{equation}
\begin{equation} \nonumber
\frac{1}{2}\sum_{i=0}^{\infty}\frac{(-1)^i\pi^{2i+2}n^{2i+1}}{2i+1}\sum_{j=0}^{i}\frac{B_{2j} n^{-2j}}{(2i+1-2j)!(2j)!}=
\frac{\pi^2}{2}\int_{0}^{n} \frac{1}{2n}\cot{\frac{\pi x}{2n}}\sin{\pi x}\,dx 
\end{equation}
\indent Let us summarize the convoluted part by summing up the three resulting integrals:
\begin{multline} \nonumber
\sum_{i=0}^{\infty}\frac{(-1)^i\pi^{2i+2}n^{2i+1}}{(2i+3)!}\sum_{j=0}^{i}\frac{(2i)! B_{2j} n^{-2j}}{(2i+1-2j)!(2j)!}=\int_{0}^{n} \frac{\pi^2(n-x)^2}{4n^3}\cot{\frac{\pi x}{2n}}\sin{\pi x}\,dx 
\\=\frac{\pi^2}{4}\int_{0}^{1}u^2\sin{\pi n(1-u)}\cot{\frac{\pi (1-u)}{2}}\,du
\end{multline}
\noindent where we made a change of variables, $u=1-x/n$.\\

\indent Now, by summing up the two parts, we obtain a formula for $H_{2}(n)$:
\begin{equation} \nonumber
\sum_{k=1}^{n}\frac{1}{k^2}=\frac{1}{2n^2}-\frac{\pi^2}{12}+\frac{\pi^2}{4}\int_{0}^{1} u^2\sin{\pi n(1-u)}\tan{\frac{\pi u}{2}}\,du \text{,} 
\end{equation}
\noindent where the identity $\cot{\pi(1-u)/2}=\tan{\pi u/2}$ was used.\\

In section \secrefe{HN_2k_sin_k_pi} we find out the general polynomial, $p_{2k}(u)$, that goes under the integral symbol, and it is convenient to move the constant $-\pi^2/12$ under the integral symbol (this also makes the $H_{2}(n)$ formula look more similar to $H(n)$).\\

First, we note that for all positive integer $n$:
\begin{equation} \nonumber
\int_{0}^{1}\sin{\pi n(1-u)}\tan{\frac{\pi u}{2}}\,du=1   \text{, which stems from the below equation:}
\end{equation}
\begin{equation} \nonumber
H^*_{0}(n)=\sum_{k=1}^{n}1=n=n+\frac{1}{2}-\frac{1}{2}\int_{0}^{1}\sin{\pi n(1-u)}\tan{\frac{\pi u}{2}}\,du=n+H_{0}(n)\text{} 
\end{equation}
\indent From the above we conclude that $H_{0}(n)=0$ for all positive integer $n$ (this is not a usual definition, but it will make sense when we reach section \secrefe{HN_2k_sin_k_pi}). Therefore our modified formula is:
\begin{equation}
\sum_{k=1}^{n}\frac{1}{k^2}=\frac{1}{2n^2}+\pi^2\int_{0}^{1}\left(-\frac{1}{12}+\frac{u^2}{4}\right)\sin{\pi n(1-u)}\tan{\frac{\pi u}{2}}\,du 
\end{equation}
\indent It can be proved that for all integer $k \geq 0$:
\begin{equation} \label{eq:integ_limit_2k} \nonumber
\lim_{n\to\infty}\int_{0}^{1}u^{2k}\sin{\pi n(1-u)}\tan{\frac{\pi u}{2}}\,du=1 \Rightarrow \lim_{n\to\infty} H_{2}(n)=\pi^2\left(-\frac{1}{12}+\frac{1}{4}\right)=\frac{\pi^2}{6}=\zeta(2) 
\end{equation}
\indent The above limits are justified by Theorem 1, section  \secrefe{Lims_1}, and Theorem 3, section \secrefe{Integ_2k}. The latter assumes that the closed-form of $\zeta(2k)$ is known, as the limits of the above integrals stem from the limits of $H_{2k}(n)$ and vice-versa. 

\subsection{Harmonic Number of Order 3} \label{HN_3_sin_k_pi}
We divide both sides of \eqrefe{start_iden_sin_k_pi_1} by $k^2$ and simplify:
\begin{equation} \nonumber
H_{3}(n)=\frac{\pi^2}{3!}H_{1}(n)-\frac{1}{2}\sum_{i=0}^{\infty}\frac{(-1)^{i}\pi^{2i+4}n^{2i+1}}{(2i+5)!}
-\sum_{i=0}^{\infty}\frac{(-1)^{i}\pi^{2i+4}n^{2i+2}}{(2i+5)!}\sum_{j=0}^{i}\frac{(2i+1)! B_{2j} n^{-2j}}{(2i+2-2j)!(2j)!} \text{,}
\end{equation}
\noindent which gives the below recurrence:
\begin{equation} \nonumber
H_{3}(n)=\frac{\pi^2}{3!}H_{1}(n)+\frac{1}{2\pi n^4}\left(\pi n-\frac{(\pi n)^3}{3!}-\sin{\pi n}\right)
-\frac{\pi^3}{12}\int_{0}^{1}u^3\left(1-\cos{\pi n(1-u)}\right)\tan{\frac{\pi u}{2}}\,dx  
\end{equation}
\indent By performing all the necessary calculations, one obtains:
\begin{equation}
\sum_{k=1}^{n}\frac{1}{k^3}=\frac{1}{2n^3}+\frac{\pi^3}{12}\int_{0}^{1}\left(u-u^3\right)\left(1-\cos{\pi n(1-u)}\right)\tan{\frac{\pi u}{2}}\,du 
\end{equation}
\indent Besides, due to the below identity, whose proof is given in section \secrefe{Lim_HN_2k+1_sin_k_pi}:
\begin{equation} \nonumber
\frac{\pi^3}{12}\int_{0}^{1}\left(u-u^3\right)\tan{\frac{\pi u}{2}}\,du=\zeta(3) \text{, the previous equation can be rewritten as:} 
\end{equation}

\begin{equation} \nonumber
\sum_{k=1}^{n}\frac{1}{k^3}=\frac{1}{2n^3}+\zeta(3)-\frac{\pi^3}{12}\int_{0}^{1} \left(u-u^3\right)\cos{\pi n(1-u)}\tan{\frac{\pi u}{2}}\,du 
\end{equation}
\indent And since the limit of $H_{3}(n)$ when $n$ tends to infinity is $\zeta(3)$, it means that:
\begin{equation} \nonumber
\lim_{n\to\infty}\int_{0}^{1}\left(u-u^3\right)\cos{\pi n(1-u)}\tan{\frac{\pi u}{2}}\,du=0 
\end{equation}

\subsection{Harmonic Number of Order 4} \label{HN_4_sin_k_pi}
We divide both sides of \eqrefe{start_iden_sin_k_pi_1} by $k^3$ and simplify:
\begin{equation} \nonumber
H_{4}(n)=\frac{\pi^2}{3!}H_{2}(n)+\frac{1}{2}\frac{\pi^4}{5!}-\frac{1}{2}\sum_{i=0}^{\infty}\frac{(-1)^i\pi^{2i+4}n^{2i}}{(2i+5)!}
-\sum_{i=0}^{\infty}\frac{(-1)^i\pi^{2i+4}n^{2i+1}}{(2i+5)!}\sum_{j=0}^{i}\frac{(2i)! B_{2j} n^{-2j}}{(2i+1-2j)!(2j)!} \text{,}
\end{equation}
\noindent which leads to the below recurrence:
\begin{equation} \nonumber
H_{4}(n)=\frac{\pi^2}{3!}H_{2}(n)+\frac{1}{2\pi n^5}\left(\pi n-\frac{(\pi n)^3}{3!}+\frac{(\pi n)^5}{5!}-\sin{\pi n}\right)-\frac{\pi^4}{48}\int_{0}^{1}u^4\sin{\pi n(1-u)}\tan{\frac{\pi u}{2}}\,du \text{} 
\end{equation}
\indent Now, after performing the calculations, we arrive at a formula for $H_{4}(n)$:
\begin{equation} \nonumber
\sum_{k=1}^{n}\frac{1}{k^4}=\frac{1}{2n^4}-\frac{7\pi^4}{720}+\pi^4\int_{0}^{1}\left(\frac{u^2}{24}-\frac{u^4}{48}\right)\sin{\pi n(1-u)}\tan{\frac{\pi u}{2}}\,du 
\end{equation}
\indent Moving the constant under the integral symbol, as we did for $H_{2}(n)$:
\begin{equation}
\sum_{k=1}^{n}\frac{1}{k^4}=\frac{1}{2n^4}+\pi^4\int_{0}^{1}\left(-\frac{7}{720}+\frac{u^2}{24}-\frac{u^4}{48}\right)\sin{\pi n(1-u)}\tan{\frac{\pi u}{2}}\,du 
\end{equation}
\indent Since each limit is 1, as mentioned in section \secrefe{HN_2_sin_k_pi}, we conclude that the limit of $H_4(n)$ is:
\begin{equation} \nonumber
\zeta(4)=\pi^4\left(-\frac{7}{720}+\frac{1}{24}-\frac{1}{48}\right)=\frac{\pi^4}{90} \text{} 
\end{equation}

\subsection{General Formula} \label{Gen_sin_k_pi}
As we have seen, lots of patterns emerge when we create formulae for $H_{1}(n)$, $H_{2}(n)$, and so on. We now assume that these patterns always repeat and see if we can find out what the general rule is for each $k$.\\

Because in each case the term that goes outside of the integral symbol is easy to deduce ($1/(2n^{2k})$ or $1/(2n^{2k+1})$), we can focus on the polynomials in $u$ that go under the integral symbol, $p_{2k}(u)$ and $p_{2k+1}(u)$, and see if we can find out their generating function, $g(x)$.\\

Note that since for each recursive equation the coefficients $\pi^{2k}$ (or $\pi^{2k+1}$) cancel out, we can ignore them for simplification purposes.

\subsection{Harmonic Numbers of Order $2k$} \label{HN_2k_sin_k_pi}
Let $f(u,n)$ be the below function (not to be confused with the same function from previous sections):
\begin{equation} \nonumber
f(u,n)=\sin{\pi n(1-u)}\tan{\frac{\pi u}{2}}
\end{equation}
\indent As we have seen in previous sections, calculating the $2k$-th harmonic number involves a recurrence with prior ones:
\begin{equation} \nonumber
H_{2k}(n)=\frac{1}{2 n^{2k}}\sum_{j=0}^{k}\frac{(-1)^{j}(\pi n)^{2j}}{(2j+1)!}-\sum_{j=0}^{k-1}\frac{(-1)^{k-j}\pi^{2k-2j}}{(2k+1-2j)!}H_{2j}(n)-\frac{(-1)^k\pi^{2k}}{2(2k)!}\int_{0}^{1}u^{2k}f(u,n)\,du \text{} 
\end{equation}
\indent That is, the harmonic numbers of even orders obey the below recursive equations:
\begingroup
\normalsize
\begin{equation} \nonumber
\begin{cases}
H_{0}(n)=\frac{1}{2}-\frac{1}{2}\int_{0}^{1}f(u,n)\,du\\

H_{2}(n)=\frac{\pi^2}{3!}H_{0}(n)+\frac{1}{2\pi n^3}\left(\pi n-\frac{\pi^3 n^3}{3!}\right)+\frac{\pi^2}{2}\int_{0}^{1}\frac{u^2}{2!}f(u,n)\,du\\

H_{4}(n)=\frac{\pi^2}{3!}H_{2}(n)-\frac{\pi^4}{5!}H_{0}(n)+\frac{1}{2\pi n^5}\left(\pi n-\frac{\pi^3 n^3}{3!}+\frac{\pi^5 n^5}{5!}\right)-\frac{\pi^4}{2}\int_{0}^{1}\frac{u^4}{4!}f(u,n)\,du\\

H_{6}(n)=\frac{\pi^2}{3!}H_{4}(n)-\frac{\pi^4}{5!}H_{2}(n)+\frac{\pi^6}{7!}H_{0}(n)+\frac{1}{2\pi n^7}\left(\pi n-\frac{\pi^3 n^3}{3!}+\frac{\pi^5 n^5}{5!}-\frac{\pi^7 n^7}{7!}\right)+\frac{\pi^6}{2}\int_{0}^{1}\frac{u^6}{6!}f(u,n)\,du\\

\vdots 
\end{cases}
\end{equation}
\endgroup
\indent Let us try to solve that recurrence, we have:
\begin{equation} \nonumber
\left(1-\frac{x^2}{3!}+\frac{x^4}{5!}-\frac{x^6}{7!}+\dotsb\right)\left(p_0+p_2 x^2+p_4 x^4+p_6 x^6+\dotsb\right)=-1+\frac{u^2}{2!}x^2-\frac{u^4}{4!}x^4-\frac{u^6}{6!}x^6+\dotsb 
\end{equation}\\
\indent But the product on the left-hand side gives:
\begin{equation} \nonumber
p_0+\left(p_2-\frac{1}{3!}p_0\right) x^2+\left(p_4-\frac{1}{3!}p_2+\frac{1}{5!}p_0\right) x^4+\left(p_6-\frac{1}{3!}p_4+\frac{1}{5!}p_2-\frac{1}{7!}p_0\right) x^6+\dotsb \text{,} 
\end{equation}
\noindent where the coefficient of each $x^{2k}$ is the recurrence that produces the polynomials $p_{2k}$ that we are interested in. The generating function for $p_{2k}(u)$ is therefore given by:
\begin{multline} \nonumber
g(x)=-\frac{x\cos{x u}}{\sin{x}}=-1+\left(-\frac{1}{6}+\frac{u^2}{2}\right)x^2+\left(-\frac{7}{360}+\frac{u^2}{12}-\frac{u^4}{24}\right)x^4\\+\left(-\frac{31}{15120}+\frac{7u^2}{720}-\frac{u^4}{144}+\frac{u^6}{720}\right)x^6+\dotsb 
\end{multline}
\indent To obtain the power series of the function $g(x)$, we need to obtain the power series of each of its components individually\citesup{Abra}: 
\begin{equation} \nonumber
\frac{x}{\sin{x}}=\sum_{i=0}^{\infty}\frac{(-1)^i B_{2i}(2-2^{2i})}{(2i)!}x^{2i} \text{, and } -{\cos{x u}}=-\sum_{i=0}^{\infty}\frac{(-1)^i u^{2i}}{(2i)!}x^{2i}
\end{equation}
\indent Therefore, the $2k$-th term of the power series of $g(x)$ is $p_{2k}(u)$, which is given by the below expression:
\begin{equation} \nonumber
p_{2k}(u)=\sum_{j=0}^{k}\frac{(-1)^{j} B_{2j}\left(2-2^{2j}\right)}{(2j)!}\cdot\frac{(-1)^{k-j+1} u^{2k-2j}}{(2k-2j)!} 
\end{equation}
\indent Now, putting it all together, we find that for all integer $k\geq 1$:
\begin{equation} \nonumber
H_{2k}(n)=\frac{1}{2n^{2k}}+\frac{\pi^{2k}}{2}\int_{0}^{1}p_{2k}(u)f(u,n)\,du \Rightarrow 
\end{equation}
\begin{equation} \label{eq:H_2k(n)}
H_{2k}(n)=\frac{1}{2n^{2k}}-\frac{(-1)^{k}\pi^{2k}}{2}\int_{0}^{1}\sum_{j=0}^{k}\frac{B_{2j}\left(2-2^{2j}\right)u^{2k-2j}}{(2j)!(2k-2j)!}\sin{\pi n(1-u)}\tan{\frac{\pi u}{2}}\,du \text{} 
\end{equation}
\indent Note that this formula applies to $H_{0}(n)$ as well, but remember that per this definition $H_{0}(n)$ is such that $H_{0}(n)=0$ for all positive integer $n$.\\

We can rewrite $H_{2k}(n)$ by means of Bernoulli polynomials\citesup{Abra}, which are given by:
\begin{equation} \nonumber
B_k(u)=\sum_{j=0}^{k}{k\choose j}B_{k-j}u^j  
\end{equation}
\indent In doing so we obtain an expression that resembles, but is not exactly, an Euler polynomial\citesup{Abra}:
\begin{equation} \nonumber
H_{2k}(n)=\frac{1}{2n^{2k}}-\frac{(-1)^{k}\pi^{2k}}{(2k)!}\int_{0}^{1}\left(B_{2k}(u)-2^{2k-1}B_{2k}\left(\frac{u}{2}\right)\right)f(u,n)\,du 
\end{equation}

\subsubsection{Generating Function of $H_{2k}(n)$}
A generating function for $H_{2k}(n)$ can be obtained by means of the generating function $g(x)$, that we previously found for $p_{2k}(u)$, as follows:
\begin{equation} \nonumber
\sum_{k=0}^{\infty}H_{2k}(n)x^{2k}=\frac{n^2}{2(n^2-x^2)}-\frac{\pi x}{2\sin{\pi x}}\int_{0}^{1}\cos{\pi x u}\sin{\pi n(1-u)}\tan{\frac{\pi u}{2}}\,du  \text{}
\end{equation}
\indent Note the convergence radius of the power series on the left-hand side is the open interval $(-1,1)$, but the domain of the function on the right-hand side is $\mathbb{R}\backslash\mathbb{Z}$. This generating function is probably an analytic continuation of the power series to the left.\\
\begin{equation} \nonumber
\indent \text{As an example, if $n=2$ the above function is } \frac{5x^2}{4}+\frac{17x^4}{16}+\frac{65x^6}{64}+\frac{257x^8}{256}+\dotsb \hspace{10cm}
\end{equation}
\indent Notice it does not have the independent term, which is a result of $H_{0}(n)=0$ for all positive integer $n$.

\subsubsection{Limit of the Generating Function of $H_{2k}(n)$}
The limit of the generating function we just found as $n$ goes to infinity is:
\begin{equation} \nonumber
h(x)=\lim_{n\to\infty}\sum_{k=1}^{\infty}H_{2k}(n)x^{2k}=\lim_{n\to\infty}\frac{n^2}{2(n^2-x^2)}-\frac{\pi x}{2\sin{\pi x}}\int_{0}^{1}\cos{\pi x u}\sin{\pi n(1-u)}\tan{\frac{\pi u}{2}}\,du \Rightarrow 
\end{equation}
\begin{equation} \nonumber
h(x)=\sum_{k=1}^{\infty}\zeta(2k)x^{2k}=\frac{1}{2}-\frac{\pi x\cos{\pi x}}{2\sin{\pi x}} \text{}
\end{equation}
\indent Note $h(x)$ also does not have the independent term, due to $H_{0}(n)=0$.\\

\textbf{Proof } The proof of the above is simple:
\begin{equation} \nonumber
\lim_{n\to\infty}\int_{0}^{1}\cos{\pi x u}\sin{\pi n(1-u)}\tan{\frac{\pi u}{2}}\,du=\lim_{n\to\infty}\int_{0}^{1}\sum_{k=0}^{\infty}\frac{(-1)^k(\pi x u)^{2k}}{(2k)!}\sin{\pi n(1-u)}\tan{\frac{\pi u}{2}}\,du \Rightarrow
\end{equation}
\begin{equation} \nonumber
\lim_{n\to\infty}\sum_{k=0}^{\infty}\frac{(-1)^k(\pi x)^{2k}}{(2k)!}\int_{0}^{1}u^{2k}\sin{\pi n(1-u)}\tan{\frac{\pi u}{2}}\,du=\cos{\pi x} \text{,} 
\end{equation}
\noindent as the limits of the above integrals are always $1$, per the two different proofs that are provided in sections \secrefe{Lims_1}, Theorem 1, and (\ref{Integ_2k}), Theorem 3. $\square$

\subsection{Harmonic Numbers of Order $2k+1$} \label{HN_2k+1_sin_k_pi}
Let $f(u,n)$ be the below function:
\begin{equation} \nonumber
f(u,n)=\left(1-\cos{\pi n(1-u)}\right)\tan{\frac{\pi u}{2}} 
\end{equation}
\indent Calculating the odd harmonic numbers involves a recurrence with prior ones:
\begingroup
\small
\begin{equation} \nonumber
H_{2k+1}(n)=\frac{1}{2 n^{2k+1}}\sum_{j=0}^{k}\frac{(-1)^{j}(\pi n)^{2j}}{(2j+1)!}-\sum_{j=0}^{k-1}\frac{(-1)^{k-j}\pi^{2k-2j}}{(2k+1-2j)!}H_{2j+1}(n)+\frac{(-1)^k\pi^{2k+1}}{2(2k+1)!}\int_{0}^{1}u^{2k+1}f(u,n)\,du 
\end{equation}
\endgroup
\indent The reasoning employed to find out the generating function of $p_{2k+1}(u)$ is entirely analogous to what we have done previously, and $g(x)$ is given by:\\
\begin{multline} \nonumber
g(x)=\frac{x\sin{x u}}{\sin{x}}=u x+\left(\frac{u}{6}-\frac{u^3}{6}\right)x^3+\left(\frac{7u}{360}-\frac{u^3}{36}+\frac{u^5}{120}\right)x^5\\+\left(\frac{31u}{15120}-\frac{7u^3}{2160}+\frac{u^5}{720}-\frac{u^7}{5040}\right)x^7+\dotsb 
\end{multline}
\indent The $(2k+1)$-th term of the power series of $g(x)$ is therefore:
\begin{equation} \nonumber
p_{2k+1}(u)=(-1)^{k}\sum_{j=0}^{k}\frac{B_{2j}\left(2-2^{2j}\right)u^{2k+1-2j}}{(2j)!(2k+1-2j)!} 
\end{equation}
\indent Now, putting it all together, we find that for all integer $k\geq 0$:
\begin{equation} \nonumber
H_{2k+1}(n)=\frac{1}{2n^{2k+1}}+\frac{\pi^{2k+1}}{2}\int_{0}^{1}p_{2k+1}(u)f(u,n)\,du \Rightarrow 
\end{equation}
\begin{equation} \nonumber
H_{2k+1}(n)=\frac{1}{2n^{2k+1}}+\frac{(-1)^{k}\pi^{2k+1}}{2}\int_{0}^{1}\sum_{j=0}^{k}\frac{B_{2j}\left(2-2^{2j}\right)u^{2k+1-2j}}{(2j)!(2k+1-2j)!}\left(1-\cos{\pi n(1-u)}\right)\tan{\frac{\pi u}{2}}\,du \text{} 
\end{equation}
\indent Because of Theorem 8, section \ref{Lim_HN_2k+1_sin_k_pi}, we can also rewrite $H_{2k+1}(n)$ as:
\begingroup
\small
\begin{equation} \nonumber
H_{2k+1}(n)=\frac{1}{2n^{2k+1}}+\zeta(2k+1)-\frac{(-1)^{k}\pi^{2k+1}}{2}\int_{0}^{1}\sum_{j=0}^{k}\frac{B_{2j}\left(2-2^{2j}\right)u^{2k+1-2j}}{(2j)!(2k+1-2j)!}\cos{\pi n(1-u)}\tan{\frac{\pi u}{2}}\,du  
\end{equation}
\endgroup
\indent Finally, with the aforementioned Bernoulli polynomials, we can also rewrite $H_{2k+1}(n)$ as:
\begin{equation} \nonumber
H_{2k+1}(n)=\frac{1}{2n^{2k+1}}+\frac{(-1)^{k}\pi^{2k+1}}{(2k+1)!}\int_{0}^{1}\left(B_{2k+1}(u)-2^{2k}B_{2k+1}\left(\frac{u}{2}\right)\right)f(u,n)\,du 
\end{equation}

\subsubsection{Generating Function of $H_{2k+1}(n)$} \label{HN_2k+1_GF}
A generating function for $H_{2k+1}(n)$ can be obtained from the function $g(x)$ that we found for $p_{2k+1}(u)$ previously, as follows:
\begin{equation} \nonumber
\sum_{k=0}^{\infty}H_{2k+1}(n)x^{2k+1}=\frac{n x}{2(n^2-x^2)}+\frac{\pi x}{2\sin{\pi x}}\int_{0}^{1}\sin{\pi x u}\left(1-\cos{\pi n(1-u)}\right)\tan{\frac{\pi u}{2}}\,du 
\end{equation}
\begin{equation} \nonumber
\indent \text{For example, if $n=2$, the above function becomes }\frac{3x}{2}+\frac{9x^3}{8}+\frac{33x^5}{32}+\frac{129x^7}{128}+\dotsb  \hspace{10cm}
\end{equation}

\subsubsection{Limit of the Generating Function of $H_{2k+1}(n)$}
\indent Before we can take the limit of this generating function as $n$ approaches infinity, we need to exclude term $H_1(n)x$, since $H(n)$ is unbounded. Hence, using the expression for $H(n)$ from section \secrefe{HN_1_sin_k_pi}, the limit of the generating function as $n$ increases is:
\begin{equation} \nonumber
\sum_{k=1}^{\infty}H_{2k+1}(n)x^{2k+1}=\frac{n x}{2(n^2-x^2)}-\frac{x}{2n}+\frac{\pi x}{2}\int_{0}^{1}\left(\frac{\sin{\pi x u}}{\sin{\pi x}}-u\right)\left(1-\cos{\pi n(1-u)}\right)\tan{\frac{\pi u}{2}}\,du \Rightarrow
\end{equation}
\begin{equation} \nonumber
h(x)=\sum_{k=1}^{\infty}\zeta(2k+1)x^{2k+1}=\frac{\pi x}{2}\int_{0}^{1}\left(\frac{\sin{\pi x u}}{\sin{\pi x}}-u\right)\tan{\frac{\pi u}{2}}\,du  \text{}
\end{equation}
\textbf{Proof } To prove that the generating function converges to the above limit, we need to show that the below integral goes to $0$ as $n$ approaches infinity. But since,
\begin{multline} \nonumber
\int_{0}^{1}\left(\sin{\pi x u}-u\sin{\pi x}\right)\cos{\pi n(1-u)}\tan{\frac{\pi u}{2}}\,du\\=\int_{0}^{1}\sum_{k=0}^{\infty}\frac{(-1)^k(\pi x)^{2k+1}\left(u^{2k+1}-u\right)}{(2k+1)!}\cos{\pi n(1-u)}\tan{\frac{\pi u}{2}}\,du \text{,}
\end{multline}
\noindent it follows that:
\begin{equation} \nonumber
\lim_{n\to\infty}\sum_{k=0}^{\infty}\frac{(-1)^k(\pi x)^{2k+1}}{(2k+1)!}\int_{0}^{1}\left(u^{2k+1}-u\right)\cos{\pi n(1-u)}\tan{\frac{\pi u}{2}}\,du=0  \text{,} 
\end{equation}
\noindent as the limits of the integrals are $0$, per Corollary 1 of section \secrefe{Lims_2}. $\square$\\

The above representation of the generating function of $\zeta(2k+1)$ is different from the one found in the literature, which employs the digamma function, though they must be equivalent:
\begin{equation} \nonumber
h(x)=\sum_{k=1}^{\infty}\zeta(2k+1)x^{2k+1}=-x\gamma-\frac{x}{2}(\psi(1+x)+\psi(1-x)) 
\end{equation}

\subsection{Initial Equation $\sin{2\pi k}=0$} \label{HN_sin_2k_pi}
In this section we set the initial equation to $\sin{2\pi k}=0$. To avoid redundancy, we omit the step by step demonstrations and present only the final formulae.\\

Using this initial equation, we obtain slightly different formulae for $H(n)$ and $H_{2}(n)$:
\begin{equation} \nonumber
\sum_{k=1}^{n}\frac{1}{k}=\frac{1}{2n}+\pi\int_{0}^{1} (1-u)\left(1-\cos{2\pi nu}\right)\cot{\pi u}\,du 
\end{equation}
\begin{equation} \nonumber
\sum_{k=1}^{n}\frac{1}{k^2}=\frac{1}{2n^2}-\frac{\pi^2}{3}-\pi^2\int_{0}^{1} u^2\sin{2\pi n(1-u)}\cot{\pi u}\,du 
\end{equation}

\subsubsection{General Formula} \label{Gen_sin_2k_pi}
We conclude that not much really changes in the system of recurrence equations, except for the introduction of a coefficient $2$ on $\pi$. Therefore, the polynomial solution is the same as before, only the multiplier of the integral and the integrand change.

\subsubsection{Harmonic Numbers of Order $2k$} \label{HN_2k_sin_2k_pi}
The recurrence equation changes slightly:
\begin{multline} \nonumber
H_{2k}(n)=\frac{1}{2 n^{2k}}\sum_{j=0}^{k}\frac{(-1)^{j}(2\pi n)^{2j}}{(2j+1)!}-\sum_{j=0}^{k-1}\frac{(-1)^{k-j}(2\pi)^{2k-2j}}{(2k+1-2j)!}H_{2j}(n)\\+\frac{(-1)^k(2\pi)^{2k}}{2(2k)!}\int_{0}^{1}u^{2k}\sin{2\pi n(1-u)}\cot{\pi u}\,du 
\end{multline}
\indent For all integer $k\geq 0$:
\begin{equation} \nonumber
H_{2k}(n)=\frac{1}{2n^{2k}}+\frac{(-1)^{k}(2\pi)^{2k}}{2}\int_{0}^{1}\sum_{j=0}^{k}\frac{B_{2j}\left(2-2^{2j}\right)u^{2k-2j}}{(2j)!(2k-2j)!}\sin{2\pi n(1-u)}\cot{\pi u}\,du \text{,} 
\end{equation}

\subsubsection{Harmonic Numbers of Order $2k+1$} \label{HN_2k+1_sin_2k_pi}
The recurrence equation also changes slightly:
\begin{multline} \nonumber
H_{2k+1}(n)=\frac{1}{2 n^{2k+1}}\sum_{j=0}^{k}\frac{(-1)^{j}(2\pi n)^{2j}}{(2j+1)!}-\sum_{j=0}^{k-1}\frac{(-1)^{k-j}(2\pi)^{2k-2j}}{(2k+1-2j)!}H_{2j+1}(n)\\-\frac{(-1)^k(2\pi)^{2k+1}}{2(2k+1)!}\int_{0}^{1}u^{2k+1}\left(1-\cos{2\pi n(1-u)}\right)\cot{\pi u}\,du 
\end{multline}
\indent For all integer $k\geq 0$:
\begin{equation} \nonumber
H_{2k+1}(n)=\frac{1}{2n^{2k+1}}-\frac{(-1)^{k}(2\pi)^{2k+1}}{2}\int_{0}^{1}\sum_{j=0}^{k}\frac{B_{2j}\left(2-2^{2j}\right)u^{2k+1-2j}}{(2j)!(2k+1-2j)!}\left(1-\cos{2\pi n(1-u)}\right)\cot{\pi u}\,du 
\end{equation}
\indent Here the integral symbol has changed due to $\cot{\pi (1-u)}=-\cot{\pi u}$.

\subsection{Initial Equation $\cos{2\pi k}=1$} \label{HN_cos_2k_pi}
When we switch to cosine-based harmonic numbers, the degree of the polynomials $p_k(u)$ go up by one. Below are two examples for $H(n)$ and $H_{2}(n)$:
\begin{equation} \nonumber
\sum_{k=1}^{n}\frac{1}{k}=\frac{1}{2n}-\pi\int_{0}^{1} u^2\left(1-\cos{2\pi n (1-u)}\right)\cot{\pi u}\,du 
\end{equation}
\begin{equation} \nonumber
\sum_{k=1}^{n}\frac{1}{k^2}=\frac{1}{2n^2}-\frac{\pi^2}{6}-\frac{2\pi^2}{3}\int_{0}^{1}u^3\sin{2\pi n (1-u)}\cot{\pi u}\,du 
\end{equation}

\subsubsection{General Formula}
Using this initial equation, we add more entropy to the formula. Here we only show a detailed demonstration for the odd case, and the even case is just stated.

\subsubsection{Harmonic Numbers of Order $2k$} \label{HN_2k_cos_2k_pi}
$H_{2k}(n)$ is given by the below recurrence equation:
\begin{multline} \nonumber
H_{2k}(n)=\frac{1}{n^{2k}}\sum_{j=0}^{k}\frac{(-1)^{j}(2\pi n)^{2j}}{(2j+2)!}-2\sum_{j=0}^{k-1}\frac{(-1)^{k-j}(2\pi)^{2k-2j}}{(2k+2-2j)!}H_{2j}(n)\\+\frac{(-1)^k(2\pi)^{2k}}{(2k+1)!}\int_{0}^{1}u^{2k+1}\sin{2\pi n(1-u)}\cot{\pi u}\,du  
\end{multline}
\indent The polynomial $p_{2k}(u)$ can be obtained using a similar approach to $p_{2k+1}(u)$ (see the next section for a detailed demonstration), which results in the formula below:
\begin{multline} \nonumber
H_{2k}(n)=\frac{1}{2n^{2k}}\\+\frac{(-1)^k \pi^{2k}}{2}\int_{0}^{1}\sum_{i=0}^{k}\sum_{j=0}^{i}\frac{B_{2j}B_{2i-2j}\left(2-2^{2j}\right)\left(2-2^{2i-2j}\right)(2u)^{2k+1-2i}}{(2j)!(2i-2j)!(2k+1-2i)!}\sin{2\pi n(1-u)}\cot{\pi u}\,du 
\end{multline}

\subsubsection{Harmonic Numbers of Order $2k+1$} \label{HN_2k+1_cos_2k_pi}
Let $f(u,n)$ be the below function:
\begin{equation} \nonumber
f(u,n)=\left(1-\cos{2\pi n(1-u)}\right)\cot{\pi u} 
\end{equation}
\indent $H_{2k+1}(n)$ is given by the below recurrence equation:
\begin{multline} \nonumber
H_{2k+1}(n)=\frac{1}{n^{2k+1}}\sum_{j=0}^{k}\frac{(-1)^{j}(2\pi n)^{2j}}{(2j+2)!}-2\sum_{j=0}^{k-1}\frac{(-1)^{k-j}(2\pi)^{2k-2j}}{(2k+2-2j)!}H_{2j+1}(n)\\-\frac{(-1)^k(2\pi)^{2k+1}}{(2k+2)!}\int_{0}^{1}u^{2k+2}f(u,n)\,du 
\end{multline}
\indent That is, the harmonic numbers of odd orders obey the below recursive equations (notice we are ignoring $\cos{2\pi n}-1$):
\begingroup
\footnotesize
\begin{equation} \nonumber
\begin{cases}
H_{1}(n)=\frac{1}{2n}-\frac{2\pi}{2!}\int_{0}^{1}u^2f(u,n)\,du\\

H_{3}(n)=2\frac{(2\pi)^2}{4!}H_{1}(n)+\frac{1}{n^3}\left(\frac{1}{2!}-\frac{(2\pi n)^2}{4!}\right)+\frac{(2\pi)^3}{4!}\int_{0}^{1}u^4f(u,n)\,du\\

H_{5}(n)=2\left(\frac{(2\pi)^2}{4!}H_{3}(n)-\frac{(2\pi)^4}{6!}H_{1}(n)\right)+\frac{1}{n^5}\left(\frac{1}{2!}-\frac{(2\pi n)^2}{4!}+\frac{(2\pi n)^4}{6!}\right)-\frac{(2\pi)^5}{6!}\int_{0}^{1}u^6 f(u,n)\,du\\

H_{7}(n)=2\left(\frac{(2\pi)^2}{4!}H_{5}(n)-\frac{(2\pi)^4}{6!}H_{3}(n)+\frac{(2\pi)^6}{8!}H_{1}(n)\right)+\frac{1}{n^7}\left(\frac{1}{2!}-\frac{(2\pi n)^4}{4!}+\frac{(2\pi n)^4}{6!}-\frac{(2\pi n)^6}{8!}\right)+\frac{(2\pi)^7}{8!}\int_{0}^{1}u^8f(u,n)\,du\\

\vdots  
\end{cases} \text{,}
\end{equation}
\endgroup
\indent Let $p(x)$ be the generating function that we are interested in. We have:
\begin{equation} \nonumber
p(x)-p(x)\frac{1}{2x^2}\left(\cos{2x}-1+\frac{(2x)^2}{2}\right)=-\sum_{k=0}^{\infty}\frac{(-1)^k 2^{2k+1}u^{2k+2}}{(2k+2)!}x^{2k+1}=\frac{1}{2x}\left(\cos{2u x}-1\right) \Rightarrow 
\end{equation}
\begin{multline} \nonumber
p_1 x+\left(-\frac{p_1}{3}+p_3\right) x^3+\left(\frac{2p_1}{45}-\frac{p_3}{3}+p_5\right) x^5+\left(-\frac{p_1}{315}+\frac{2p_3}{45}-\frac{p_5}{3}+p_7\right) x^7+\dotsb=\\
-u^2 x+\frac{u^4}{3}x^3-\frac{2u^6}{45}x^5+\frac{u^8}{315}x^7+\dotsb  
\end{multline}
\noindent The generating function for $p_{2k+1}(u)$ is therefore given by:
\begin{multline} \nonumber
g(x)=\left(\frac{x}{\sin{x}}\right)^2\frac{\cos{2u x}-1}{2x}=-u^2 x+\left(-\frac{u^2}{3}+\frac{u^4}{3}\right)x^3+\left(-\frac{u^2}{15}+\frac{u^4}{9}-\frac{2u^6}{45}\right)x^5\\+\left(-\frac{2u^2}{189}+\frac{u^4}{45}-\frac{2u^6}{135}+\frac{u^8}{315}\right)x^7+\dotsb \text{,} 
\end{multline}
\noindent where we used following transformation:
\begin{equation} \nonumber
\noindent -x\frac{\cos{2u x}-1}{\cos{2x}-1+(2x)^2/2-2x^2}=-x\frac{\cos{2x u}-1}{\cos{2x}-1}=x\,\frac{\cos{2u x}-1}{2\sin^2{x}}
\end{equation}
\indent To obtain the power series of $g(x)$, we need to obtain the power series of each function individually:
\begin{equation} \nonumber
\left(\frac{x}{\sin{x}}\right)^2=\sum_{i=0}^{\infty}\sum_{j=0}^{i}\frac{(-1)^i B_{2j}B_{2i-2j}\left(2-2^{2j}\right)\left(2-2^{2i-2j}\right)}{(2j)!(2i-2j)!}x^{2i}  \text{, } \frac{\cos{2x u}-1}{x}=\sum_{i=1}^{\infty}\frac{(-1)^i (2u)^{2i}}{(2i)!}x^{2i-1}
\end{equation}
\indent Therefore, the $(2k+1)$-th term of the power series of $g(x)$ is given by:
\begin{equation} \nonumber
p_{2k+1}(u)=\frac{1}{2}\sum_{i=0}^{k}\sum_{j=0}^{i}\frac{(-1)^i B_{2j}B_{2i-2j}\left(2-2^{2j}\right)\left(2-2^{2i-2j}\right)}{(2j)!(2i-2j)!}\cdot\frac{(-1)^{k+1-i} (2u)^{2k+2-2i}}{(2k+2-2i)!} \text{,} 
\end{equation}
\noindent which goes into the final formula:
\begin{equation} \nonumber
H_{2k+1}(n)=\frac{1}{2n^{2k+1}}+\pi^{2k+1}\int_{0}^{1}p_{2k+1}(u)\left(1-\cos{2\pi n(1-u)}\right)\cot{\pi u}\,du 
\end{equation}

\section{Alternating Harmonic Numbers: $C^2_{k}(n)$} \label{HN_cos_k_pi}
Setting the initial equation to $\cos{\pi k}=(-1)^k $  drastically changes the picture. It no longer enables to calculate $H_{k}(n)$, but the alternating harmonic numbers instead, $C^2_{k}(n)$.\\

Below are a couple of examples of formulae for the alternating harmonic numbers:
\begin{equation} \nonumber
C^2_{1}(n)=\sum_{k=1}^{n}\frac{(-1)^k}{k}=H_{1}(n)+\frac{1}{2n}\left(-1+\cos{\pi n}\right)-\frac{\pi}{2}\int_{0}^{1}\left(1-\cos{\pi n(1-u)}\right)\tan{\frac{\pi u}{2}}\,du 
\end{equation}
\begin{equation} \nonumber
C^2_{2}(n)=\sum_{k=1}^{n}\frac{(-1)^k}{k^2}=H_{2}(n)+\frac{1}{2n^2}\left(-1+\frac{(\pi n)^2}{2!}+\cos{\pi n}\right)-\frac{\pi^2}{2}\int_{0}^{1}u\sin{\pi n(1-u)}\tan{\frac{\pi u}{2}}\,du 
\end{equation}

\subsection{General Formula: $C^2_{k}(n)$} \label{Gen_cos_k_pi}

The recurrence equations for the generalized alternating harmonic numbers are:
\begin{multline} \nonumber
C^2_{2k}(n)=\sum_{j=1}^{n}\frac{(-1)^j}{j^{2k}}=\frac{1}{2n^{2k}}\left(\cos{\pi n}-\sum_{j=0}^{k}\frac{(-1)^j (\pi n)^{2j}}{(2j)!}\right)\\+\sum_{j=0}^{k}\frac{(-1)^{k-j}\pi^{2k-2j}}{(2k-2j)!}H_{2j}(n)+\frac{(-1)^k\pi^{2k}}{2(2k-1)!}\int_{0}^{1} u^{2k-1}\sin{\pi n(1-u)}\tan{\frac{\pi u}{2}}\,du 
\end{multline}
\begin{multline} \nonumber
C^2_{2k+1}(n)=\sum_{j=1}^{n}\frac{(-1)^j}{j^{2k+1}}=\frac{1}{2n^{2k+1}}\left(\cos{\pi n}-\sum_{j=0}^{k}\frac{(-1)^j (\pi n)^{2j}}{(2j)!}\right)\\+\sum_{j=0}^{k}\frac{(-1)^{k-j}\pi^{2k-2j}}{(2k-2j)!}H_{2j+1}(n)-\frac{(-1)^k\pi^{2k+1}}{2(2k)!}\int_{0}^{1} u^{2k}\left(1-\cos{\pi n(1-u)}\right)\tan{\frac{\pi u}{2}}\,du 
\end{multline}

\section{Odd Alternating Harmonic Numbers: $S^4_{k}(n)$} \label{HN_sin_k/2_pi}
If we set the initial equation to $\sin{\pi k/2}$, we are able to obtain formulae for the odd alternating harmonic numbers, $S^4_{k}(n)$.\\

Two examples of formulae for the odd alternating harmonic numbers are below:
\begin{equation} \nonumber
S^4_{1}(n)=\sum_{k=1}^{n}\frac{1}{k}sin{\frac{\pi k}{2}}=\frac{1}{2n}\left(-\frac{\pi n}{2}+\sin{\frac{\pi n}{2}}\right)+\frac{\pi}{4}\int_{0}^{1}\sin{\frac{\pi n(1-u)}{2}}\left(\sec{\frac{\pi u}{2}}+\tan{\frac{\pi u}{2}}\right)\,du 
\end{equation}
\begin{multline} \nonumber
S^4_{2}(n)=\sum_{k=1}^{n}\frac{1}{k^2}\sin{\frac{\pi k}{2}}=\frac{\pi}{2}H_{1}(n)+\frac{1}{2 n^2}\left(-\frac{\pi n}{2}+\sin{\frac{\pi n}{2}}\right)\\-\frac{\pi^2}{8}\int_{0}^{1} u\left(1-\cos{\frac{\pi n(1-u)}{2}}\right)\left(\sec{\frac{\pi u}{2}}+\tan{\frac{\pi u}{2}}\right)\,du 
\end{multline}

\subsection{General Formula: $S^4_{k}(n)$} \label{Gen_sin_k/2_pi}
The recurrence equations for the generalized odd alternating harmonic numbers are:
\begin{multline} \nonumber
S^4_{2k}(n)=\sum_{j=1}^{n}\frac{1}{j^{2k}}\sin{\frac{\pi j}{2}}=\frac{1}{2n^{2k}}\left(\sin{\frac{\pi n}{2}}-\sum_{j=0}^{k-1}\frac{(-1)^j (\frac{\pi n}{2})^{2j+1}}{(2j+1)!}\right)\\-\sum_{j=0}^{k-1}\frac{(-1)^{k-j}(\frac{\pi}{2})^{2k-1-2j}}{(2k-1-2j)!}H_{2j+1}(n)\\+\frac{(-1)^k(\frac{\pi}{2})^{2k}}{2(2k-1)!}\int_{0}^{1} u^{2k-1}\left(1-\cos{\frac{\pi n(1-u)}{2}}\right)\left(\sec{\frac{\pi u}{2}}+\tan{\frac{\pi u}{2}}\right)\,du  
\end{multline}
\begin{multline} \nonumber
S^4_{2k+1}(n)=\sum_{j=1}^{n}\frac{1}{j^{2k+1}}\sin{\frac{\pi j}{2}}=\frac{1}{2n^{2k+1}}\left(\sin{\frac{\pi n}{2}}-\sum_{j=0}^{k}\frac{(-1)^j (\frac{\pi n}{2})^{2j+1}}{(2j+1)!}\right)\\+\sum_{j=0}^{k}\frac{(-1)^{k-j}(\frac{\pi}{2})^{2k+1-2j}}{(2k+1-2j)!}H_{2j}(n)\\
+\frac{(-1)^k(\frac{\pi}{2})^{2k+1}}{2(2k)!}\int_{0}^{1} u^{2k}\sin{\frac{\pi n(1-u)}{2}}\left(\sec{\frac{\pi u}{2}}+\tan{\frac{\pi u}{2}}\right)\,du  
\end{multline}

\section{General Formula: $C^z_{k}(n)$ and $S^z_{k}(n)$}
There is a striking similarity between the formulae derived from initial equation $\cos{\pi k}=1$ and the ones derived with $\sin{\pi k/2}$. Based on this similarity, we are able to generalize the pattern.

\subsection{$C^z_{2k}(n)$ and $S^z_{2k+1}(n)$} \label{Final_1}
We grouped these two under the same section because they share an integral and they both have $H_{2j}(n)$ in their recursions.\\

For all complex $z$, $C^z_{2k}(n)$ and $S^z_{2k+1}(n)$ are given by:
\begin{multline} \label{eq:C^z_(2k)(n)}
C^z_{2k}(n)=\sum_{j=1}^{n}\frac{1}{j^{2k}}\cos{\frac{2\pi j}{z}}=\frac{1}{2n^{2k}}\left(\cos{\frac{2\pi n}{z}}-\sum_{j=0}^{k}\frac{(-1)^j (\frac{2\pi n}{z})^{2j}}{(2j)!}\right)\\+\sum_{j=0}^{k}\frac{(-1)^{k-j}(\frac{2\pi}{z})^{2k-2j}}{(2k-2j)!}H_{2j}(n)
\\+\frac{(-1)^k(\frac{2\pi}{z})^{2k}}{2(2k-1)!}\int_{0}^{1}(1-u)^{2k-1}\sin{\frac{2\pi n u}{z}}\cot{\frac{\pi u}{z}}\,du \text{, }\forall\text{ integer }k \geq 1
\end{multline}

\begin{multline} \label{eq:S^z_(2k+1)(n)}
S^z_{2k+1}(n)=\sum_{j=1}^{n}\frac{1}{j^{2k+1}}\sin{\frac{2\pi j}{z}}=\frac{1}{2n^{2k+1}}\left(\sin{\frac{2\pi n}{z}}-\sum_{j=0}^{k}\frac{(-1)^j (\frac{2\pi n}{z})^{2j+1}}{(2j+1)!}\right)\\+\sum_{j=0}^{k}\frac{(-1)^{k-j}(\frac{2\pi}{z})^{2k+1-2j}}{(2k+1-2j)!}H_{2j}(n)\\+\frac{(-1)^k(\frac{2\pi}{z})^{2k+1}}{2(2k)!}\int_{0}^{1}(1-u)^{2k}\sin{\frac{2\pi n u}{z}}\cot{\frac{\pi u}{z}}\,du \text{, }\forall\text{ integer }k \geq 0
\end{multline}
\indent Notice that in order for these equations to hold, we need to have $H_{0}(n)=0$ for all positive integer $n$, as mentioned before.

\subsubsection{Limits of $C^z_{2k}(n)$ and $S^z_{2k+1}(n)$} \label{Lims_1}

At infinity, $C^z_{2k}(n)$ and $S^z_{2k+1}(n)$ become Fourier series, denoted here by $C^z_{2k}$ and $S^z_{2k+1}$, whose closed-forms are given by Bernoulli polynomials, per Abramowitz and Stegun:\citesup{Abra}
\begin{equation} \nonumber
\sum_{j=1}^{\infty}\frac{1}{j^{2k}}\cos{2\pi x j}=\frac{-(-1)^k(2\pi)^{2k}}{2(2k)!}B_{2k}(x) \text{ and } \sum_{j=1}^{\infty}\frac{1}{j^{2k+1}}\sin{2\pi x j}=\frac{(-1)^k(2\pi)^{2k+1}}{2(2k+1)!}B_{2k+1}(x)
\end{equation}
\indent The above result implies the following theorem, which holds for all integer $k \geq 0$ and real $z\geq 1$:
\begin{equation} \nonumber
\textbf{Theorem 1}\lim_{n\to\infty}\int_{0}^{1}(1-u)^{k}\sin{\frac{2\pi n u}{z}}\cot{\frac{\pi u}{z}}\,du=
\begin{cases}
      1, & \text{if}\ k=0 \text{ and }z=1\\
      z/2, & \text{otherwise}
\end{cases} 
\end{equation}

Therefore, with the exception of $S^1_{1}=0$, the limits of $C^z_{2k}(n)$ and $S^z_{2k+1}(n)$, for real $z \geq 1$, are given by:
\begin{equation} \nonumber
C^z_{2k}=\sum_{j=1}^{\infty}\frac{1}{j^{2k}}\cos{\frac{2\pi j}{z}}=\sum_{j=0}^{k}\frac{(-1)^{k-j}\left(\frac{2\pi}{z}\right)^{2k-2j}}{(2k-2j)!}\zeta(2j)+\frac{(-1)^k z}{4(2k-1)!}\left(\frac{2\pi}{z}\right)^{2k} \text{ (} \forall \text{ }k\geq 1\text{)} 
\end{equation}
\begin{equation} \nonumber
S^z_{2k+1}=\sum_{j=1}^{\infty}\frac{1}{j^{2k+1}}\sin{\frac{2\pi j}{z}}=\sum_{j=0}^{k}\frac{(-1)^{k-j}\left(\frac{2\pi}{z}\right)^{2k+1-2j}}{(2k+1-2j)!}\zeta(2j)+\frac{(-1)^k z}{4(2k)!}\left(\frac{2\pi}{z}\right)^{2k+1} \text{ (} \forall\text{ } k\geq 0\text{)} 
\end{equation}
\indent These are just rewrites of the expressions for  $C^z_{2k}$ and $S^z_{2k+1}$ from reference \citena{Abra}, with $x=1/z$.

\subsection{$C^z_{2k+1}(n)$ and $S^z_{2k}(n)$} \label{Final_2}
For all complex $z$, $C^z_{2k+1}(n)$ and $S^z_{2k}(n)$ are given by:
\begin{multline} \label{eq:C^z_(2k+1)(n)}
C^z_{2k+1}(n)=\sum_{j=1}^{n}\frac{1}{j^{2k+1}}\cos{\frac{2\pi j}{z}}=\frac{1}{2n^{2k+1}}\left(\cos{\frac{2\pi n}{z}}-\sum_{j=0}^{k}\frac{(-1)^j (\frac{2\pi n}{z})^{2j}}{(2j)!}\right)
\\+\sum_{j=0}^{k}\frac{(-1)^{k-j}(\frac{2\pi}{z})^{2k-2j}}{(2k-2j)!}H_{2j+1}(n)\\-\frac{(-1)^k(\frac{2\pi}{z})^{2k+1}}{2(2k)!}\int_{0}^{1}(1-u)^{2k}\left(1-\cos{\frac{2\pi n u}{z}}\right)\cot{\frac{\pi u}{z}}\,du \text{, }\forall\text{ integer }k \geq 0
\end{multline}

\begin{multline} \label{eq:S^z_(2k)(n)}
S^z_{2k}(n)=\sum_{j=1}^{n}\frac{1}{j^{2k}}\sin{\frac{2\pi j}{z}}=\frac{1}{2n^{2k}}\left(\sin{\frac{2\pi n}{z}}-\sum_{j=0}^{k-1}\frac{(-1)^j (\frac{2\pi n}{z})^{2j+1}}{(2j+1)!}\right)
\\-\sum_{j=0}^{k-1}\frac{(-1)^{k-j}(\frac{2\pi}{z})^{2k-1-2j}}{(2k-1-2j)!}H_{2j+1}(n)\\+\frac{(-1)^k(\frac{2\pi}{z})^{2k}}{2(2k-1)!}\int_{0}^{1}(1-u)^{2k-1}\left(1-\cos{\frac{2\pi n u}{z}}\right)\cot{\frac{\pi u}{z}}\,du \text{, }\forall\text{ integer }k \geq 1
\end{multline}

\subsubsection{Limits of $C^z_{2k+1}(n)$ and $S^z_{2k}(n)$} \label{Lims_2}

Before taking the limit of $C^z_{2k+1}(n)$ and $S^z_{2k}(n)$ as $n$ tends to infinity, we need to remove $H(n)$ from the second sum (since it explodes out to infinity), and recombine it with the integral.\\

\indent In order to do that, we need to use one of the three formulae we created for $H(n)$ in sections \secrefe{HN_1_sin_k_pi}, \secrefe{HN_sin_2k_pi} and \secrefe{HN_cos_2k_pi}. Since the last two are almost identical, let us consider only the first two:
\begin{equation} \nonumber
H(n)-\frac{1}{2n}=\frac{\pi}{2}\int_{0}^{1}(1-u)\left(1-\cos{\pi n u}\right)\cot{\frac{\pi u}{2}}\,du=\pi\int_{0}^{1}(1-u)\left(1-\cos{2\pi n u}\right)\cot{\pi u}\,du 
\end{equation}
\indent By using either one of these formulae, we can carve out two integrals, one that does not depend on $n$ and one that does, as in the below example:
\begin{multline} \nonumber
C^z_{2k+1}(n)=\frac{1}{2n^{2k+1}}\left(\cos{\frac{2\pi n}{z}}-\sum_{j=0}^{k}\frac{(-1)^j (\frac{2\pi n}{z})^{2j}}{(2j)!}\right)+\sum_{j=1}^{k}\frac{(-1)^{k-j}(\frac{2\pi}{z})^{2k-2j}}{(2k-2j)!}H_{2j+1}(n)\\-\frac{(-1)^k(\frac{2\pi}{z})^{2k+1}}{2(2k)!}
\bigg(-\frac{z}{2\pi n}+\int_{0}^{1}(1-u)^{2k}\cot{\frac{\pi u}{z}}-z(1-u)\cot{\pi u}\,du \\-\int_{0}^{1}(1-u)^{2k}\cos{\frac{2\pi nu}{z}}\cot{\frac{\pi u}{z}}-z (1-u)\cos{2\pi n u}\cot{\pi u}\,du \bigg)
\end{multline}
\indent Therefore, to know the limit of $C^z_{2k+1}(n)$, we need to know the limit of the integral to the right as $n$ grows. This limit is given in the following theorem, which holds for all integer $k \geq 0$ and real $z\geq 1$ (except $k=0$ and $z=1$), and for which we do not provide a proof:
\begin{equation} \nonumber
\textbf{Theorem 2}
\lim_{n\to\infty}\int_{0}^{1}(1-u)^{k}\cos{\frac{2\pi n u}{z}}\cot{\frac{\pi u}{z}}-z(1-u)\cos{2\pi n u}\cot{\pi u}\,du=\frac{z\log(z)}{\pi} 
\end{equation}
\indent This limit apparently does not exist in the literature.
Theorem 2 allows to deduce the following corollary:
\begin{equation} \nonumber
\textbf{Corollary 1}
\lim_{n\to\infty}\int_{0}^{1}\left(u^{k}-u\right)\cos{\pi n (1-u)}\tan{\frac{\pi u}{2}}\,du=0 \text{ } \forall\text{ integer }k\geq 0 
\end{equation}
\noindent \textbf{Proof 1} This result stems from Theorem 2 and the fact we can write $C^z_{2k+1}(n)$ or $S^z_{2k}(n)$ using different formulae for $H(n)$, which leads to an equation:
\begin{multline} \nonumber
\int_{0}^{1}(1-u)^{k}\cos{\frac{2\pi n u}{z}}\cot{\frac{\pi u}{z}}-\frac{z}{2}(1-u)\cos{\pi n u}\cot{\frac{\pi u}{2}}\,du-\int_{0}^{1}(1-u)^k\cot{\frac{\pi u}{z}}-\frac{z}{2}(1-u)\cot{\frac{\pi u}{2}}\,du
\\=\int_{0}^{1}(1-u)^{k}\cos{\frac{2\pi n u}{z}}\cot{\frac{\pi u}{z}}-z (1-u)\cos{2\pi n u}\cot{\pi u}\,du-\int_{0}^{1}(1-u)^k\cot{\frac{\pi u}{z}}-z(1-u)\cot{\pi u}\,du 
\end{multline}
\begin{multline} \nonumber
\Rightarrow \int_{0}^{1}(1-u)^{k}\cos{\frac{2\pi n u}{z}}\cot{\frac{\pi u}{z}}-\frac{z}{2}(1-u)\cos{\pi n u}\cot{\frac{\pi u}{2}}\,du\\=\int_{0}^{1}(1-u)^{k}\cos{\frac{2\pi n u}{z}}\cot{\frac{\pi u}{z}}-z (1-u)\cos{2\pi n u}\cot{\pi u}\,du+\int_{0}^{1}z(1-u)\cot{\pi u}-\frac{z}{2}(1-u)\cot{\frac{\pi u}{2}}\,du
\end{multline}
\indent Now, by making $z=2$ and using Theorem 2, it follows from the above relation that:
\begin{equation} \nonumber
\lim_{n\to\infty}\int_{0}^{1}\left(u^{k}-u\right)\cos{\pi n(1-u)}\tan{\frac{\pi u}{2}}\,du=\frac{2\log(2)}{\pi}+\int_{0}^{1}2(1-u)\cot{\pi u}-(1-u)\cot{\frac{\pi u}{2}}\,du=0\text{. } \square
\end{equation}\\
\indent Now that we have Theorem 2, we can find out the limits of $C^z_{2k+1}(n)$ and $S^z_{2k}(n)$ which, except for $C^1_1=\infty$, are given by:
\begingroup
\small
\begin{multline} \label{eq:C^z_(2k+1)}
C^z_{2k+1}=\sum_{j=1}^{\infty}\frac{1}{j^{2k+1}}\cos{\frac{2\pi j}{z}}=\sum_{j=1}^{k}\frac{(-1)^{k-j}\left(\frac{2\pi}{z}\right)^{2k-2j}}{(2k-2j)!}\zeta(2j+1)+\frac{(-1)^k \log(z)\left(\frac{2\pi}{z}\right)^{2k}}{(2k)!}\\-\frac{(-1)^k(\frac{2\pi}{z})^{2k+1}}{2(2k)!}\int_{0}^{1}(1-u)^{2k}\cot{\frac{\pi u}{z}}-z(1-u)\cot{\pi u}\,du 
\end{multline}
\begin{multline} \label{eq:S^z_(2k)}
S^z_{2k}=\sum_{j=1}^{\infty}\frac{1}{j^{2k}}\sin{\frac{2\pi j}{z}}=-\sum_{j=1}^{k-1}\frac{(-1)^{k-j}\left(\frac{2\pi}{z}\right)^{2k-1-2j}}{(2k-1-2j)!}\zeta(2j+1)-\frac{(-1)^k \log(z)\left(\frac{2\pi}{z}\right)^{2k-1}}{(2k-1)!}\\
+\frac{(-1)^k(\frac{2\pi}{z})^{2k}}{2(2k-1)!}\int_{0}^{1}(1-u)^{2k-1}\cot{\frac{\pi u}{z}}-z(1-u)\cot{\pi u}\,du
\end{multline}
\endgroup
\indent Note that $H(n)$ diverges because $\int_{0}^{1}\cot{\pi u}-(1-u)\cot{\pi u}\,du$ diverges.

\section{Example: Infinite Sum of $H(n)/n^2$}

In this section, we derive expressions for sums of the type $H_{k}(n)/n^{r}$, over the positive integers $n$, with $k$ odd and $r$ even, and vice-versa. We do not try to derive the result for $k$ and $r$ both even or odd, because these cases lead to integrals that are very hard to evaluate.\\ 

Hence, let us start with an example. We want to obtain the sum of $H(n)/n^{2}$ over the positive integers using the formula for $H(n)$ from section \secrefe{HN_sin_2k_pi}:\\
\begin{equation} \nonumber
H(n)=\frac{1}{2n}+\pi\int_{0}^{1} u\left(1-\cos{2\pi n(1-u)}\right)\cot{\pi(1-u)}\,du \Rightarrow  
\end{equation}
\begin{multline} \nonumber
\sum_{n=1}^{\infty}\frac{H(n)}{n^2}=\sum_{n=1}^{\infty}\frac{1}{n^2}\left(\frac{1}{2n}+\pi\int_{0}^{1} (1-u)\left(1-\cos{2\pi n u}\right)\cot{\pi u}\,du \right)=\\
\frac{1}{2}\zeta(3)+\pi\int_{0}^{1} (1-u)\left(\zeta(2)-\sum_{n=1}^{\infty}\frac{1}{n^2}\cos{2\pi n u}\right)\cot{\pi u}\,du
\end{multline}
\indent The Fourier series can be simplified using the results from section \secrefe{Lims_1}, giving:
\begin{equation} \nonumber
\sum_{n=1}^{\infty}\frac{H(n)}{n^2}=\frac{1}{2}\zeta(3)+\pi^3\int_{0}^{1} u(1-u)^2\cot{\pi u}\,du=2\,\zeta(3)  
\end{equation}

\subsection{General Formula: Sum of $H_{2k}(n)/n^{2r+1}$}

Here we use the formula for $H_{2k}(n)$ from section (\ref{HN_2k_sin_2k_pi}), with a slight transformation only valid for integer $n$:
\begin{equation} \nonumber
\sum_{n=1}^{\infty}\frac{H_{2k}(n)}{n^{2r+1}}=\sum_{n=1}^{\infty}\frac{1}{n^{2r+1}}\left(\frac{1}{2n^{2k}}-\frac{(-1)^{k}(2\pi)^{2k}}{2}\int_{0}^{1}\sum_{j=0}^{k}\frac{B_{2j}\left(2-2^{2j}\right)u^{2k-2j}}{(2j)!(2k-2j)!}\sin{2\pi n u}\cot{\pi u}\,du\right) 
\end{equation}
\begin{equation} \nonumber
\sum_{n=1}^{\infty}\frac{H_{2k}(n)}{n^{2r+1}}=\frac{\zeta(2k+2r+1)}{2}-\frac{(-1)^{k}(2\pi)^{2k}}{2}\int_{0}^{1}\sum_{j=0}^{k}\frac{B_{2j}\left(2-2^{2j}\right)u^{2k-2j}}{(2j)!(2k-2j)!}\left(\sum_{n=1}^{\infty}\frac{\sin{2\pi n u}}{n^{2r+1}}\right)\cot{\pi u}\,du 
\end{equation}
\indent Now, the closed-form of the Fourier series above, after a change of variables, $u=1/z$, is given by:
\begin{equation} \nonumber
\sum_{n=1}^{\infty}\frac{\sin{2\pi n u}}{n^{2r+1}}=\sum_{i=0}^{r}\frac{(-1)^{r-i}\zeta(2i)(2\pi u)^{2r+1-2i}}{(2r+1-2i)!}+\frac{(-1)^r (2\pi)^{2r+1}u^{2r}}{4(2r)!}=-\frac{(-1)^{r}(2\pi)^{2r+1}}{2(2r+1)!}B_{2r+1}(u) 
\end{equation}
\indent Therefore, we can express the sum as function of Bernoulli polynomials, and it is finite for all integer $r\geq 1$:
\begin{multline} \nonumber
\sum_{n=1}^{\infty}\frac{H_{2k}(n)}{n^{2r+1}}=\frac{\zeta(2k+2r+1)}{2}\\+\frac{(-1)^{k+r}(2\pi)^{2k+2r+1}}{2(2k)!(2r+1)!}\int_{0}^{1}\left(B_{2k}(u)-2^{2k-1}B_{2k}\left(\frac{u}{2}\right)\right)B_{2r+1}(u)\cot{\pi u}\,du 
\end{multline}
\indent If $k=0$ the sum is always zero (since $H_{0}(n)=0$), which enables to deduce another integral representation for $\zeta(2r+1)$, which happens to coincide with the one in Abramowitz and Stegun\citesup{Abra}:
\begin{equation} \nonumber
\zeta(2r+1)=-\frac{(-1)^{r}(2\pi)^{2r+1}}{2(2r+1)!}\int_{0}^{1}B_{2r+1}(u)\cot{\pi u}\,du 
\end{equation}

\subsection{General Formula: Sum of $H_{2k+1}(n)/n^{2r}$}

Here we use the formula for $H_{2k+1}(n)$ from section \secrefe{HN_2k+1_sin_2k_pi}, though we could have used the two others as well (notice we made a transformation only valid for integer $n$):
\begin{multline} \nonumber
\sum_{n=1}^{\infty}\frac{H_{2k+1}(n)}{n^{2r}}=\sum_{n=1}^{\infty}\frac{1}{n^{2r}}\bigg(\frac{1}{2n^{2k+1}}\\-\frac{(-1)^{k}(2\pi)^{2k+1}}{2}\int_{0}^{1}\sum_{j=0}^{k}\frac{B_{2j}\left(2-2^{2j}\right)u^{2k+1-2j}}{(2j)!(2k+1-2j)!}\left(1-\cos{2\pi n u}\right)\cot{\pi u}\,du\bigg)
\end{multline}
\begin{multline} \nonumber
\sum_{n=1}^{\infty}\frac{H_{2k+1}(n)}{n^{2r}}=\frac{\zeta(2k+1+2r)}{2}\\-\frac{(-1)^{k}(2\pi)^{2k+1}}{2}\int_{0}^{1}\sum_{j=0}^{k}\frac{B_{2j}\left(2-2^{2j}\right)u^{2k+1-2j}}{(2j)!(2k+1-2j)!}\left(\sum_{n=1}^{\infty}\frac{1-\cos{2\pi n u}}{n^{2r}}\right)\cot{\pi u}\,du 
\end{multline}
\indent Now, the closed-form of the Fourier series featured in the above equation is given in section \secrefe{Lims_1}, and it can also be expressed as Bernoulli polynomials. That is, after a change of variables, $u=1/z$, we obtain:
\begin{equation} \nonumber
\sum_{n=1}^{\infty}\frac{\cos{2\pi n u}}{n^{2r}}=\sum_{i=0}^{r}\frac{(-1)^{r-i}\zeta(2i)(2\pi u)^{2r-2i}}{(2r-2i)!}+\frac{(-1)^r (2\pi)^{2r}u^{2r-1}}{4(2r-1)!}=-\frac{(-1)^{r}(2\pi)^{2r}}{2(2r)!}B_{2r}(u) 
\end{equation}
\indent In a way, the closed-form from section \secrefe{Lims_1} is more general than the Bernoulli polynomial. For instance, when the denominator here is $n^{2r+1}$ instead of $n^{2r}$, we can use its analogous form from section \secrefe{Lims_2}, which is no longer a Bernoulli polynomial.\\

\indent That said, for integer $k\geq 0$ and $r \geq 1$, we can write:
\begin{multline}  \nonumber
\sum_{n=1}^{\infty}\frac{H_{2k+1}(n)}{n^{2r}}=\frac{\zeta(2k+1+2r)}{2}\\
-\frac{(-1)^{k}(2\pi)^{2k+1}}{(2k+1)!}\int_{0}^{1}\left(B_{2k+1}(u)-2^{2k}B_{2k+1}\left(\frac{u}{2}\right)\right)\left(\zeta(2r)+\frac{(-1)^{r}(2\pi)^{2r}}{2(2r)!}B_{2r}(u)\right)\cot{\pi u}\,du 
\end{multline}
\indent Or, although we lose the validity of the formula for $k=0$, we can rewrite the sum as:
\begin{multline} \nonumber
\sum_{n=1}^{\infty}\frac{H_{2k+1}(n)}{n^{2r}}=\frac{\zeta(2k+1+2r)}{2}+\zeta(2k+1)\zeta(2r)\\-\frac{(-1)^{k+r}(2\pi)^{2k+1+2r}}{2(2k+1)!(2r)!}\int_{0}^{1}\left(B_{2k+1}(u)-2^{2k}B_{2k+1}\left(\frac{u}{2}\right)\right)B_{2r}(u)\cot{\pi u}\,du 
\end{multline}

\section{Limits of the Integrals} 

\subsection{Limits of the Integrals in the $H_{2k}(n)$ Recursions} \label{Integ_2k}

In this section we present proofs for the limits of the integrals that appear in the recursions of $H_{2k}(n)$ from sections \secrefe{HN_2k_sin_k_pi}, \secrefe{HN_2k_sin_2k_pi} and \secrefe{HN_2k_cos_2k_pi}. This approach requires prior knowledge of the closed-forms of $\zeta(2k)$, as mentioned in sections \secrefe{HN_2_sin_k_pi} and \secrefe{HN_4_sin_k_pi}.\\

Looking back at the set of recurrence equations from the aforementioned sections, it is evident that we can express each integral as a function of $H_{2j}(n)$:
\begingroup
\small
\begin{equation} \nonumber
\int_{0}^{1}u^{2k}\sin{\pi n(1-u)}\tan{\frac{\pi u}{2}}\,du=\frac{-2(-1)^k(2k)!}{\pi^{2k}}\left(\sum_{j=1}^{k}\frac{(\pi\ii)^{2k-2j}}{(2k+1-2j)!}H_{2j}(n)-\frac{1}{2 n^{2k}}\sum_{j=0}^{k}\frac{(\pi\ii n)^{2j}}{(2j+1)!}\right) \text{}
\end{equation}
\begin{equation} \nonumber
\int_{0}^{1}u^{2k}\sin{2\pi n(1-u)}\cot{\pi u}\,du=\frac{2(-1)^k(2k)!}{(2\pi)^{2k}}\left(\sum_{j=1}^{k}\frac{(2\pi\ii)^{2k-2j}}{(2k+1-2j)!}H_{2j}(n)-\frac{1}{2 n^{2k}}\sum_{j=0}^{k}\frac{(2\pi\ii n)^{2j}}{(2j+1)!}\right) 
\end{equation}
\begin{equation} \nonumber
\int_{0}^{1}u^{2k+1}\sin{2\pi n(1-u)}\cot{\pi u}\,du=\frac{2(-1)^k(2k+1)!}{(2\pi)^{2k}}\left(\sum_{j=1}^{k}\frac{(2\pi\ii)^{2k-2j}}{(2k+2-2j)!}H_{2j}(n)-\frac{1}{2 n^{2k}}\sum_{j=0}^{k}\frac{(2\pi\ii n)^{2j}}{(2j+2)!}\right) 
\end{equation}
\endgroup
\indent We can deduce the limits of these integrals based on the closed-forms of $\zeta(2k)$. Conversely, knowing these limits allows to deduce the closed-forms of $\zeta(2k)$.
\begin{equation} \nonumber
\textbf{Theorem 3}\lim_{n\to\infty}\int_{0}^{1}u^{2k}\sin{\pi n(1-u)}\tan{\frac{\pi u}{2}}\,du=1 \text{ } \forall \text{ integer }k\geq 0 
\end{equation}\\
\noindent $\textbf{Proof 3}$ This integral appears with initial equation $\sin{\pi k}=0$ and per Theorem 1, section \secrefe{Lims_1}, its limit is $z/2=1$, which we shall confirm now. By taking the limit of the integral as $n$ approaches infinity, we have:
\begin{multline} \nonumber
\lim_{n\to\infty}\int_{0}^{1}u^{2k}\sin{\pi n(1-u)}\tan{\frac{\pi u}{2}}\,du=\frac{-2(-1)^k(2k)!}{\pi^{2k}}\left(\sum_{j=1}^{k}\frac{(\pi\ii)^{2k-2j}}{(2k-2j+1)!}\zeta(2j)-\frac{(\pi\ii)^{2k}}{2(2k+1)!}\right) \\=\frac{-2(-1)^k(2k)!}{\pi^{2k}}\sum_{j=0}^{k}\frac{(-1)^{k-j}\pi^{2k-2j}}{(2k-2j+1)!}\zeta(2j)=(2k)!\sum_{j=0}^{k}\frac{2^{2j}B_{2j}}{(2k-2j+1)!(2j)!}  
\end{multline}
\indent (Note that $H_{0}(n)=0$, but $\zeta(0)=-1/2$.) Now, to complete the proof, let us show that the above sum equals $1$ for all integer $k\geq 0$. For that, let $g(x)$ be the product of the two below functions:
\begin{equation} \nonumber
x\coth{x}=x\,\frac{e^{x}+e^{-x}}{e^{x}-e^{-x}}=\sum_{j=0}^{\infty}\frac{2^{2j}B_{2j}}{(2j)!}x^{2j} \text{, and} \sinh{x}=\frac{e^{x}-e^{-x}}{2}=\sum_{j=0}^{\infty}\frac{1}{(2j+1)!}x^{2j+1} \Rightarrow
\end{equation}
\begin{equation} \nonumber
g(x)=x\coth{x}\,\sinh{x}=\sum_{k=0}^{\infty}\sum_{j=0}^{k}\frac{2^{2j}B_{2j}}{(2j)!}x^{2j}\cdot \frac{1}{(2k-2j+1)!}x^{2k-2j+1} \Rightarrow  
\end{equation}
\begin{equation} \nonumber
g(x)=\sum_{k=0}^{\infty}\left(\sum_{j=0}^{k}\frac{2^{2j}B_{2j}}{(2k-2j+1)!(2j)!}\right)x^{2k+1}=x\,\frac{e^{x}+e^{-x}}{2}=x\cosh{x}=\sum_{k=0}^{\infty}\frac{1}{(2k)!}x^{2k+1} \text{,} 
\end{equation}
\noindent which implies the theorem. $\square$

\begin{equation} \nonumber
\textbf{Theorem 4}\lim_{n\to\infty}\int_{0}^{1}u^{k}\sin{2\pi n(1-u)}\cot{\pi u}\,du=
\begin{cases}
      -1, & \text{if}\ k=0 \\
      -1/2, & \text{if integer }k\geq 1
\end{cases} 
\end{equation}
\noindent $\textbf{Proof 4}$ This integral appears with initial equations $\sin{2\pi k}=0$ and $\cos{2\pi k}=1$. Here we only prove the case $\sin{2\pi k}=0$, though case $\cos{2\pi k}=1$ should follow a similar reasoning and be straightforward.\\

By taking the limit of the integral as $n$ goes to infinity, we have:
\begin{multline} \nonumber
\lim_{n\to\infty}\int_{0}^{1}u^{2k}\sin{2\pi n(1-u)}\cot{\pi u}\,du=\frac{2(-1)^k(2k)!}{(2\pi)^{2k}}\left(\sum_{j=1}^{k}\frac{(2\pi\ii)^{2k-2j}}{(2k-2j+1)!}\zeta(2j)-\frac{(2\pi\ii)^{2k}}{2(2k+1)!}\right)
\\=\frac{2(-1)^k(2k)!}{(2\pi)^{2k}}\sum_{j=0}^{k}\frac{(-1)^{k-j}(2\pi)^{2k-2j}}{(2k-2j+1)!}\zeta(2j)=-(2k)!\sum_{j=0}^{k}\frac{B_{2j}}{(2k-2j+1)!(2j)!}
\end{multline}
\indent Now, to complete the proof, let us show that the above sum equals $-1$ if $k=0$, or $-1/2$ if $k\geq 1$. For that, let $g(x)$ be the product of the two below functions, that also appeared in the proof of Theorem 3, only now the cotangent is re-scaled:
\begin{equation} \nonumber
\frac{x}{2}\coth{\frac{x}{2}}=\sum_{j=0}^{\infty}\frac{B_{2j}}{(2j)!}x^{2j} \text{, and} \sinh{x}=\sum_{j=0}^{\infty}\frac{1}{(2j+1)!}x^{2j+1}
\end{equation}
\indent Therefore, we have:
\begin{equation} \nonumber
g(x)=\frac{x}{2}\coth{\frac{x}{2}}\,\sinh{x}=\sum_{k=0}^{\infty}\sum_{j=0}^{k}\frac{B_{2j}}{(2j)!}x^{2j}\cdot \frac{1}{(2k-2j+1)!}x^{2k-2j+1} \Rightarrow 
\end{equation}
\begin{equation} \nonumber
g(x)=\sum_{k=0}^{\infty}\left(\sum_{j=0}^{k}\frac{B_{2j}}{(2k-2j+1)!(2j)!}\right)x^{2k+1}=\frac{x}{2}\left(1+\frac{e^{x}+e^{-x}}{2}\right)=x+\sum_{k=1}^{\infty}\frac{1}{2(2k)!}x^{2k+1} \text{,} 
\end{equation}
\noindent which implies the theorem. $\square$

\subsubsection{Limit of $H_{2k}(n)$} \label{Lim_HN_2k}

The limit of $H_{2k}(n)$ as $n$ approaches infinity is $\zeta(2k)$ for all integer $k\geq 0$, which is proved in the two following theorems:
\begin{equation} \nonumber
\textbf{Theorem 5}\lim_{n\to\infty}\frac{1}{2n^{2k}}-\frac{(-1)^{k}\pi^{2k}}{2}\int_{0}^{1}\sum_{j=0}^{k}\frac{B_{2j}\left(2-2^{2j}\right)u^{2k-2j}}{(2j)!(2k-2j)!}\sin{\pi n(1-u)}\tan{\frac{\pi u}{2}}\,du=\zeta(2k) 
\end{equation}
\noindent \textbf{Proof 5} First we note that due to Theorem 3, the above limit reduces to: 
\begin{equation} \nonumber
-\frac{(-1)^{k}\pi^{2k}}{2}\sum_{j=0}^{k}\frac{B_{2j}\left(2-2^{2j}\right)}{(2j)!(2k-2j)!}
\end{equation}
\indent Now we can prove that the above expression equals $\zeta(2k)$, using the same approach from the previous section:
\begin{multline} \nonumber
\sum_{k=0}^{\infty}x^{2k}\sum_{j=0}^{k}\frac{B_{2j}\left(2-2^{2j}\right)}{(2j)!(2k-2j)!}=\sum_{k=0}^{\infty}\sum_{j=0}^{k}\frac{\left(2-2^{2j}\right)B_{2j}x^{2j}}{(2j)!}\cdot\frac{x^{2k-2j}}{(2k-2j)!}\\=\left(x\coth{\frac{x}{2}}-x\coth{x}\right)\cosh{x}
=x\coth{x}=\sum_{k=0}^{\infty}\frac{2^{2k}B_{2k}}{(2k)!}x^{2k} \Rightarrow \\-\frac{(-1)^{k}\pi^{2k}}{2}\sum_{j=0}^{k}\frac{B_{2j}\left(2-2^{2j}\right)}{(2j)!(2k-2j)!}=-\frac{(-1)^{k}(2\pi)^{2k}B_{2k}}{2(2k)!}=\zeta(2k)\text{ } \square
\end{multline}
\begin{equation} \nonumber
\textbf{Theorem 6}
\lim_{n\to\infty}\frac{1}{2n^{2k}}+\frac{(-1)^{k}(2\pi)^{2k}}{2}\int_{0}^{1}\sum_{j=0}^{k}\frac{B_{2j}\left(2-2^{2j}\right)u^{2k-2j}}{(2j)!(2k-2j)!}\sin{2\pi n(1-u)}\cot{\pi u}\,du=\zeta(2k) 
\end{equation}
\noindent \textbf{Proof 6} : First we note that due to Theorem 4, the above limit reduces to:
\begin{equation} \nonumber
-\frac{(-1)^{k}(2\pi)^{2k}}{2}\left(\frac{B_{2k}(2-2^{2k})}{(2k)!}+\frac{1}{2}\sum_{j=0}^{k-1}\frac{B_{2j}\left(2-2^{2j}\right)}{(2j)!(2k-2j)!}\right)
\end{equation}
\indent Now let us prove that the above expression equals $\zeta(2k)$, using some of the previous results:
\begin{multline} \nonumber
\frac{B_{2k}(2-2^{2k})}{(2k)!}+\frac{1}{2}\sum_{j=0}^{k-1}\frac{B_{2j}\left(2-2^{2j}\right)}{(2j)!(2k-2j)!}=\frac{1}{2}\frac{B_{2k}(2-2^{2k})}{(2k)!}+\frac{1}{2}\sum_{j=0}^{k}\frac{B_{2j}\left(2-2^{2j}\right)}{(2j)!(2k-2j)!} 
\\=\frac{1}{2}\frac{B_{2k}(2-2^{2k})}{(2k)!}+\frac{1}{2}\frac{B_{2k}2^{2k}}{(2k)!}=\frac{B_{2k}}{(2k)!} \Rightarrow -\frac{(-1)^{k}(2\pi)^{2k}}{2}\,\frac{B_{2k}}{(2k)!}=\zeta(2k)\text{ }\square
\end{multline}

\subsection{Limits of the Integrals in the $H_{2k+1}(n)$ Recursions} \label{Integ_2k+1}

We can express each integral as a function of $H_{2j+1}(n)$ for all real $n$:
\begingroup
\footnotesize
\begin{equation} \nonumber
\int_{0}^{1}u^{2k+1}\left(1-\cos{\pi n(1-u)}\right)\tan{\frac{\pi u}{2}}\,du=\frac{2(-1)^k(2k+1)!}{\pi^{2k+1}}\left(\sum_{j=0}^{k}\frac{(\pi\ii)^{2k-2j}}{(2k+1-2j)!}H_{2j+1}(n)-\frac{1}{2 n^{2k+1}}\sum_{j=0}^{k}\frac{(\pi\ii n)^{2j}}{(2j+1)!}\right) 
\end{equation}
\begin{equation} \nonumber
\int_{0}^{1}u^{2k+1}\left(1-\cos{2\pi n(1-u)}\right)\cot{\pi u}\,du=\frac{-2(-1)^k(2k+1)!}{(2\pi)^{2k+1}}\left(\sum_{j=0}^{k}\frac{(2\pi\ii)^{2k-2j}}{(2k+1-2j)!}H_{2j+1}(n)-\frac{1}{2 n^{2k+1}}\sum_{j=0}^{k}\frac{(2\pi\ii n)^{2j}}{(2j+1)!}\right) 
\end{equation}
\begin{equation} \nonumber
\int_{0}^{1}u^{2k+2}\left(1-\cos{2\pi n(1-u)}\right)\cot{\pi u}\,du=\frac{-2(-1)^k(2k+2)!}{(2\pi)^{2k+1}}\left(\sum_{j=0}^{k}\frac{(2\pi\ii)^{2k-2j}}{(2k+2-2j)!}H_{2j+1}(n)-\frac{1}{2 n^{2k+1}}\sum_{j=0}^{k}\frac{(2\pi\ii n)^{2j}}{(2j+2)!}\right) 
\end{equation}
\endgroup\\
\indent From the above equations, we can infer that each one of the integrals tends to plus or minus infinity as $n$ increases (the reason is that each one contains $H(n)$, which is unbounded).\\

\indent One consequence of this fact is that the coefficients of $p_{2k+1}(u)$ in the formulae of $H_{2k+1}(n)$ need to sum up to $0$ for all $k\geq 1$, in order to cancel out those infinities, the exception being $H(n)$. This statement is translated in the next theorem.
\begin{equation} \nonumber
\textbf{Theorem 7 }p_{2k+1}(1)=\sum_{j=0}^{k}\frac{B_{2j}\left(2-2^{2j}\right)}{(2j)!(2k+1-2j)!}=0 \text{, }\forall\text{ integer } k\geq 1 
\end{equation}
\noindent \textbf{Proof 7} In section \secrefe{HN_2k+1_sin_k_pi}, we have created a generating function for $p_{2k+1}(u)$, which allows to deduce the following equivalence (notice the second sum in the double-sum is $p_{2k+1}(u)$):
\begin{equation} \nonumber
\sum_{k=0}^{\infty}(-1)^{k}x^{2k+1}\sum_{j=0}^{k}\frac{B_{2j}\left(2-2^{2j}\right)u^{2k+1-2j}}{(2j)!(2k+1-2j)!}=\frac{x\sin{x u}}{\sin{x}}  
\end{equation}
\indent Now, if $u=1$, we conclude that the above sum equals $x$, which implies the theorem. $\square$\\ 

Similarly, the below also holds, though the proof is omitted:
\begin{equation} \nonumber
\sum_{i=0}^{k}\sum_{j=0}^{i}\frac{B_{2j}B_{2i-2j}\left(2-2^{2j}\right)\left(2-2^{2i-2j}\right)2^{2k+2-2i}}{(2j)!(2i-2j)!(2k+2-2i)!}=0 \text{ } \forall \text{ }k\geq 1 
\end{equation}

\subsubsection{Limit of $H_{2k+1}(n)$} \label{Lim_HN_2k+1_sin_k_pi}

The values of $\zeta(2k+1)$ are given by the first part of the integral, as explained by the next theorem:
\begin{equation} \nonumber
\textbf{Theorem 8 }\zeta(2k+1)=\frac{(-1)^{k}\pi^{2k+1}}{2}\int_{0}^{1}\sum_{j=0}^{k}\frac{B_{2j}\left(2-2^{2j}\right)u^{2k+1-2j}}{(2j)!(2k+1-2j)!}\tan{\frac{\pi u}{2}}\,du \text{, }\forall\text{ integer } k\geq 1\text{}  
\end{equation}
\textbf{Proof 8} To prove this result, we need to show that:
\begin{equation} \nonumber
\lim_{n\to\infty}\int_{0}^{1}\sum_{j=0}^{k}\frac{B_{2j}\left(2-2^{2j}\right)u^{2k+1-2j}}{(2j)!(2k+1-2j)!}\cos{\pi n(1-u)}\tan{\frac{\pi u}{2}}\,du=0 
\end{equation}
\indent Using the result from Theorem 2, section \secrefe{Lims_2}, we know that for large $n$ we can write:
\begin{equation} \nonumber
\int_{0}^{1}u^{2k+1-2j}\cos{\pi n(1-u)}\tan{\frac{\pi u}{2}}+2u\cos{2\pi n(1-u)}\cot{\pi u}\,du \sim \frac{2\log(2)}{\pi} 
\end{equation}
\indent But per Theorem 7:
\begin{equation} \nonumber
\frac{2\log(2)}{\pi}\sum_{j=0}^{k}\frac{B_{2j}\left(2-2^{2j}\right)}{(2j)!(2k+1-2j)!}=0 \text{ and } \int_{0}^{1}\sum_{j=0}^{k}\frac{B_{2j}\left(2-2^{2j}\right)u}{(2j)!(2k+1-2j)!}\cos{2\pi n(1-u)}\cot{\pi u}\,du=0  
\end{equation}
\noindent which implies the theorem. $\square$\\

There is a slightly different integral representation for $\zeta(2k+1)$, which stems from the formula derived in section \secrefe{HN_sin_2k_pi}:
\begin{equation} \nonumber
\zeta(2k+1)=-\frac{(-1)^{k}(2\pi)^{2k+1}}{2}\int_{0}^{1}\sum_{j=0}^{k}\frac{B_{2j}\left(2-2^{2j}\right)u^{2k+1-2j}}{(2j)!(2k+1-2j)!}\cot{\pi u}\,du 
\end{equation}

\end{document}